\documentclass[a4paper,english,11pt]{amsart}

\usepackage[utf8]{inputenc}
\usepackage[T1]{fontenc}
\usepackage{babel}
\usepackage[babel]{csquotes}

\usepackage{amsmath}
\usepackage{amsfonts}
\usepackage{amsthm}
\usepackage{bbm}
\usepackage{amssymb}
\usepackage[top=3.4cm, bottom=3.6cm, left=3.6 cm, right=3.6 cm]{geometry}
\usepackage{mathrsfs}
\usepackage[pdftex, pdfstartview=FitW,colorlinks=true,linkcolor=blue,%
urlcolor=blue,citecolor=blue]{hyperref}
\usepackage{stmaryrd}
\usepackage{cleveref}
\usepackage{mathtools}
\usepackage{soul}
\usepackage{bbm}
\usepackage{bm}
\usepackage{verbatim}
\usepackage{url}
\usepackage{xcolor}
\usepackage{enumitem}
\usepackage{stmaryrd}
\usepackage{lmodern}

\newtheorem{th.}{Theorem}[section]
\newtheorem{thm}[th.]{Theorem}

\newtheorem{lemma}[th.]{Lemma} 
 
\newtheorem{proposition}[th.]{Proposition}

\newtheorem{thmx}{Theorem}

\theoremstyle{definition}
\newtheorem{definition}[th.]{Definition}
\newtheorem*{remark}{Remark}

\numberwithin{equation}{section}

\newcommand{\Q}{\mathbb{Q}}
\newcommand{\N}{\mathbb{N}}
\newcommand{\Z}{\mathbb{Z}}
\newcommand{\C}{\mathbb{C}}
\newcommand{\R}{\mathbb{R}}
\newcommand{\G}{\Gamma}

\newcommand{\kt}{\mathfrak{t}}
\newcommand{\kg}{\mathfrak{g}}

\newcommand{\kp}{\mathfrak{p}}
\newcommand{\ku}{\mathfrak{u}}
\newcommand{\ks}{\mathfrak{s}}

\newcommand {\cB} {{\mathcal B}}
\newcommand {\cC} {{\mathcal C}}
\newcommand {\cD} {{\mathcal D}}
\newcommand {\cE} {{\mathcal E}}
\newcommand {\cF} {{\mathcal F}}

\newcommand {\cL} {{\mathcal L}}

\newcommand {\cN} {{\mathcal N}}
\newcommand {\cO} {{\mathcal O}}
\newcommand {\cP} {{\mathcal P}}

\newcommand {\cS} {{\mathcal S}}

\newcommand {\cV} {{\mathcal V}}
\newcommand {\cZ} {{\mathcal Z}}

\DeclareMathOperator{\supp}{supp}

\DeclareMathOperator{\SL}{SL}

\DeclareMathOperator{\GL}{GL}

\DeclareMathOperator{\SO}{SO}

\DeclareMathOperator{\Ad}{Ad}

\newcommand{\bA}{{\bf A}}
\newcommand{\bB}{{\bf B}}
\newcommand{\bD}{{\bf D}}
\newcommand{\bG}{{\bf G}}
\newcommand{\bL}{{\bf L}}
\newcommand{\bM}{{\bf M}}
\newcommand{\bS}{{\bf S}}
\newcommand{\bT}{{\bf T}}
\newcommand{\bU}{{\bf U}}

\newcommand{\bP}{{\bf P}}
\newcommand{\bV}{{\bf V}}
\newcommand{\bX}{{\bf X}}

\newcommand{\dd}{\,\mathrm{d}}

\DeclareMathOperator{\vol}{vol}

\renewcommand{\setminus}{\smallsetminus}
\renewcommand{\subset}{\subseteq}

\title[Counting rational approximations on flag varieties]{Counting rational approximations on rank one flag varieties}

\author{Ren\'e Pfitscher}

\address{Ren\'e Pfitscher. Department of Mathematics, Universit\'e Sorbonne Paris Nord}

\email{pfitscher@math.univ-paris13.fr}

\thanks{}

\begin{document}

\begin{abstract}
On the set of real points of a generalized flag variety of rank one, we count rational approximations to almost every point with respect to the Riemannian measure. In particular, our results apply to Grassmann varieties and projective quadrics. The proof uses exponential mixing in the space of lattices and tools from geometry of numbers.
\end{abstract}

\maketitle
\setcounter{tocdepth}{1}

\section{Introduction}
Let $\psi : \N \rightarrow (0,+\infty)$ be a non-increasing function. Khintchine's theorem~\cite{Khintchine26} asserts that the inequality
\[
0 \leq q x - p < \psi(q)
\]
admits infinitely (resp. finitely) many solutions $(p,q) \in \Z \times \N$ for almost every $x \in \R$, if the series $\sum_{q = 1}^{\infty} \psi(q)$ diverges (resp. converges). In the case the series is divergent, Schmidt~\cite{Schmidt60} strengthened Khintchine's theorem. More precisely, for every $x \in \R$ and $T \geq 1$, he considered the counting function
\begin{equation} \label{eq:SchmidtCountingFunction}
\cN_\psi(x,T) = \# \left \{ (p,q) \in \Z \times \N : 0 \leq q x - p < \psi(q), \, 1 \leq q < T \right \}
\end{equation}
and showed that for almost every $x \in \R$, $\cN_\psi(x,T)$ is asymptotically equal to $\sum_{1\leq q < T} \psi(q)$ as $T$ goes to infinity, with an explicit error term. In fact, Schmidt's result holds not only for the real line, but also for the Euclidean space $\R^n$ of any dimension $n \geq 1$. 

\medskip
Our goal is to prove a version of this theorem, where the Euclidean space $\R^n$ is replaced by the set of real points $X = \bX(\R)$ of a generalized flag variety $\bX$ defined over $\Q$. A height function $H : \bX(\Q) \rightarrow (0,+\infty)$ on the set of rational points takes the place of the denominator $q$ of a rational number $\tfrac{p}{q}$ written in reduced form (see Section \ref{sec:Height}). We measure the distance between a real point and its rational approximations using a Riemannian distance $d(\cdot,\cdot)$ on $X$. 

\medskip
Let $c > 0$ and $\tau \in [0,+\infty)$ be arbitrary. In analogy to \eqref{eq:SchmidtCountingFunction}, for every $x \in X$ and $T \geq 1$, we define
\begin{equation} \label{eq:CountingFunction}
\cN_{c,\tau}(x,T) = \# \{ v \in \bX(\Q) : d(x,v) < c \, H(v)^{-\tau}, \, 1 \leq H(v) < T\}.
\end{equation}
In this paper, we determine the asymptotic behavior of $\cN_{c,\tau}(x,T)$ as $T \rightarrow + \infty$, for almost every $x \in X$ with respect to the Riemannian measure on $X$.

\subsection{Main results}
The natural setting to state our main theorem is that of generalized flag varieties, that is, varieties that can be written as the quotient $\bX = \bG / \bP$ of a connected semisimple algebraic $\Q$-group $\bG$ by a parabolic $\Q$-subgroup $\bP$ of $\bG$. Throughout this paper, we assume that $\bP$ is a maximal parabolic $\Q$-subgroup of $\bG$ with abelian unipotent radical. In particular, $\bX$ has $\Q$-rank $1$. We may, without loss of generality, assume that the group $\bG$ is simply connected and almost $\Q$-simple (see Section \ref{sec:Flag-Variety}). We denote algebraic varieties defined over $\mathbb{Q}$ by bold letters and their sets of real points by ordinary letters. Furthermore, for simplicity of exposition, we will work with the set of complex points of an algebraic variety defined over $\Q$, and refer to it simply as the variety itself when no confusion arises. For instance, we write $G = \bG(\R)$ and $\bG = \bG(\C)$ to denote the groups of real and complex points of $\bG$, respectively.

\medskip
To state our main result, we need to introduce some more notation. Let $K$ be a maximal compact subgroup of $G$. Let $\sigma_X$ be the unique $K$-invariant probability measure on $X$. We equip $X$ with a $K$-invariant Riemannian distance $d(\cdot, \cdot)$ and the set of rational points $\bX(\Q)$ with a height function $H_\chi$ associated to an irreducible rational representation $\pi_{\chi} : \bG \rightarrow \GL(\bV_{\chi})$ which is generated by a unique rational line $\bD_{\chi}$ of highest weight $\chi$ such that $\mathrm{Stab}_{\bG} \, \bD_{\chi} = \bP$ (see Section \ref{sec:Height}). By \cite[Th\'eor\`emes 2.4.5 et 3.2.1]{deSaxce20}, there exists a rational number $\beta_\chi > 0$ such that, for every $c > 0$ and for almost every $x \in X$, the inequality
\begin{equation} \label{eq:DiophantineExpo}
d(x,v) < c \, H_\chi (v)^{-\tau}
\end{equation}
admits infinitely (resp. finitely) many solutions $v \in \bX(\Q)$, if $\tau \leq \beta_\chi$ (resp. $\tau > \beta_\chi$). We refer to $\beta_\chi$ as the \emph{Diophantine exponent} of $X$ with respect to $\chi$. 

\medskip
In this paper, we strengthen this statement by providing an asymptotic formula, with an explicit error term if $\tau \in (0, \beta_{\chi})$, for the number of solutions to \eqref{eq:DiophantineExpo} with height at most $T$. Let $d = \dim X$ be the dimension of $X$ as a real manifold. For every function $\psi : \R_{+}^{\times} \rightarrow \R_{+}^{\times}$, element $x \in X$, and parameter $T \geq 1$, we define 
\[
\cN_\psi(x,T) = \# \left \{ v \in \bX(\Q) : d(x,v) < \psi(H(v)), \, 1 \leq H(v) < T \right \}.
\]

\begin{thmx} [Counting below the Diophantine exponent] \label{thm:ThmMain2}
Let $\bG$ be a connected simply connected almost $\Q$-simple $\Q$-group, $\bP$ be a maximal parabolic $\Q$-subgroup of $\bG$ with abelian unipotent radical and $\bX = \bG / \bP$. Let $H_{\chi}$ be a height function on $\bX(\Q)$ associated to a dominant $\Q$-weight $\chi$. Let $\psi : \R_{+}^{\times} \rightarrow \R_{+}^{\times}$ be a non-increasing function satisfying that there exist $\tau \in (0,\beta_\chi)$ and $C > 1$ such that $C^{-1} y^{-\tau} \leq \psi(y) \leq C y^{-\tau}$ for all sufficiently large $y \geq 1$. For every $T \geq 1$, define the integral
\[
\Psi(T) = \int_{1}^T \psi(y)^{d} y^{\beta_\chi d} \frac{d y}{y}.
\]
Then there exists an explicit constant $\varkappa > 0$ and $\varepsilon > 0$ such that, for almost every $x \in X$, as $T \rightarrow + \infty$,
\begin{equation} \label{eq:ThmMain2}
\cN_{\psi}(x, T) = \varkappa \, \Psi(T) \left ( 1 + O_x(\Psi(T)^{- \varepsilon}) \right ).
\end{equation}
\end{thmx}

\begin{remark}
The constant $\varkappa$ in Theorem \ref{thm:ThmMain2} is explicitly given in Equation \eqref{eq:Varkappa-Constant}. In particular, it depends on $\beta_{\chi}$, $d$ and the total Riemannian volume of $X$.
\end{remark}

Throughout this paper, the notation $f(T) \sim g(T)$ for real-valued functions $f$ and $g$ means that $\frac{f(T)}{g(T)} \rightarrow 1$ as $T \rightarrow + \infty$.

\begin{thmx} [Counting at the Diophantine exponent] \label{thm:ThmMain}
Let $\bG$, $\bP$, $\bX$, $\chi$, and $\varkappa > 0$ be as in Theorem \ref{thm:ThmMain2}. Let $c > 0$ and let $\tau = \beta_{\chi}$. 
Then, for almost every $x \in X$, as $T \rightarrow + \infty$,
\begin{equation} \label{eq:ThmMain}
\cN_{c,\beta_{\chi}}(x, T) \sim \varkappa \, c^d \, \ln (T).
\end{equation}
\end{thmx}

\begin{remark}
Our method of proof differs depending on whether $\tau < \beta_\chi$ or $\tau = \beta_\chi$. To handle the first case, our proof was inspired by the work of Huang and de Saxc\'e~\cite{dSH25} on the local distribution of rational points on flag varieties. For the second case, we generalize the ergodic-theoretic approach of Alam and Ghosh~\cite{AG22}, who counted rational approximations on spheres, which in turn builds on ideas from Athreya-Parrish-Tseng~\cite{APT16}. 
\end{remark}

Finally, we show that rational points are equidistributed in $X$. In \cite[Theorem~4]{MSG14}, Mohammadi and Salehi Golsefidy counted rational points of bounded height on a flag variety $\bX = \bG / \bP$ with respect to an arbitrary metrized line bundle and arbitrary parabolic $\Q$-subgroup $\bP$. We complement this statement by showing that rational points are in fact effectively equidistributed in $X$. More precisely, for every $x \in X$ and $r > 0$, we denote by $B_X(x,r)$ the open ball in $X$ with center $x$ and radius $r$, and define the counting function 
\begin{align}\label{equ:khincountfun}
\cN_{\chi}(x,r,T) =\#\left\{v \in \bX(\Q) : v \in B_X(x,r), \, 1 \leq H_\chi(v) < T\right\}.
\end{align}

\begin{thmx} [Effective equidistribution of rational points]\label{thm:ThmMainEquid}
Let $\bG$, $\bP$, $\bX$ and $\chi$ be as in Theorem \ref{thm:ThmMain2}. Then, there exists an explicit constant $\varkappa'>0$ (see Equation \ref{eq:Varkappa-Constant-Equid}) and $\varepsilon' > 0$ such that, for all $r > 0$ and $x \in X$, as $T \rightarrow + \infty$,
\begin{equation} \label{eq:ThmMainEquid}
\cN_{\chi}(x, r, T) = \varkappa' \, T^{\beta_\chi d} \sigma_X(B_X(x,r))\left ( 1 + O_r(T^{- \varepsilon'}) \right ).
\end{equation}
\end{thmx}

\begin{remark}
The proof of Theorem \ref{thm:ThmMainEquid} follows the line of argument of Theorem \ref{thm:ThmMain2}, but requires fewer ingredients. In fact, we note that for every $r > 0$, $x \in X$ and $T \geq 1$, using the notation of Theorem \ref{thm:ThmMain2} with $\psi_r : \R_{+}^{\times} \to \R_{+}^{\times}$ the constant function $\psi_r(y) = r$, we have
\[
\cN_{\chi}(x, r, T) = \cN_{\psi_r}(x, T).
\]
\end{remark}

\subsection{Context and applications}
All of the above results apply to Grassmann varieties and projective quadric hypersurfaces. After providing some context, we showcase several applications of Theorems \ref{thm:ThmMain2} and \ref{thm:ThmMain}.

\medskip
In the special case where $X$ is a sphere or a rational ellipsoid, this problem has received much attention recently. In particular, Alam and Ghosh \cite{AG22} and Ouaggag \cite{Ouaggag23} counted rational approximations on the sphere of arbitrary dimension $n$, but only for approximation functions $\psi: \N \to \R_+^\times$ of the form $\psi(q) = c q^{-1}$ with $c > 0$. 
On the other hand, Kelmer and Yu \cite{KY23} were able to deal with a general approximation function, but their result does not apply to spheres of dimension $n > 1$ with $n \equiv 1 \mod 8$. 

\medskip
More generally, let $n \in \N^\times$ and let $X_Q$ be the real points of an $n$-dimensional projective rational quadric hypersurface $\bX_Q$, given as the set of zeros in $\mathbb{P}(\R^{n+2})$ of a non-degenerate rational quadratic form $Q$ in ${n+2}$ variables:
\begin{equation} \label{eq:Def-Quadric}
X_Q = [Q^{-1}(0)] = \big\{ x \in \mathbb{P}(\R^{n+2}) : \exists \, \bm{x} \in \mathbb{R}^{n+2} \smallsetminus \{\bm{0}\}, \, x = [\bm{x}]\text{ and }Q(\bm{x}) = 0 \big\}.
\end{equation}
Here $[\bm{x}]$ denotes the projective point corresponding to a non-zero vector $\bm{x} \in \mathbb{R}^{n+2}$. The distance $d(\cdot,\cdot)$ and the height function $H$ are obtained by restriction of the usual distance and height function on $\mathbb{P}(\R^{n+2})$, respectively. Let $K$ be a maximal compact subgroup of the special orthogonal group $\SO_Q(\R)$ associated to $Q$ and let $\sigma_{Q}$ be the $K$-invariant probability measure on $X_Q$. Furthermore, one assumes that $X_Q$ contains a rational point; by stereographic projection, this implies in fact that $\bX_Q(\Q)$ is dense in $X_Q$.
Fishman, Kleinbock, Merrill, and Simmons \cite{FKMS22} have obtained the first remarkable results for intrinsic Diophantine approximation in this setting. 

\medskip
As a corollary of Theorems \ref{thm:ThmMain2} and \ref{thm:ThmMain}, we obtain the following statement, which extends the above-mentioned results. Contrary to the Kelmer-Yu theorem, we need no congruence condition on the dimension $n$, and we may in fact take $X_Q$ to be any non-degenerate projective quadric hypersurface as considered by Fishman-Kleinbock-Merrill-Simmons. Moreover, let $X_0 \subseteq \mathbb{P}(\R^4)$ be the quadric hypersurface associated to the quadratic form $Q_0(\bm{x}) = x_0x_3 - x_1x_2$. The $\Q$-rank of $X_0$ is $2$ and the formulation of the Khintchine-type theorem for quadric hypersurfaces \cite[Theorem~6.3]{FKMS22} depends on whether or not $Q$ is conjugate over $\Q$ to $Q_0$. Using \cite[Theorem~1.5]{FKMS22}, for every $c > 0$ and for almost every $x \in X_Q$, the inequality $d(x,v) < c \, H(v)^{-\tau}$ has infinitely (resp. finitely) many solutions $v \in \bX_Q(\Q)$, if $\tau \leq 1$ (resp. $\tau > 1$). 

\begin{thm} \label{thm:Quadrics}
Let $n\in \N^\times$ and let $Q$ be a non-degenerate $\Q$-isotropic rational quadratic form in $n+2$ variables that is not conjugate over $\Q$ to $Q_0$. Let $\tau \in (0,1]$ and $c > 0$. Then there exists an explicit constant $\varkappa > 0$ such that for $\sigma_{Q}$-almost every $x \in X_Q$, as $T \rightarrow + \infty$,
\begin{equation*} \label{eq:Cor1_ThmMain}
\cN_{c,\tau}(x, T) \sim \varkappa \, c^n \, \int_{1}^T y^{(1-\tau)n} \, \tfrac{d y}{y}.
\end{equation*}
\end{thm}

\begin{remark} \label{remark:Exceptional}
Unfortunately, it seems that our method does not yield a similar result for the exceptional quadric hypersurface $X_0$. We will provide more details in Section \ref{sec:Quadric}, after introducing the necessary notation throughout the paper.
\end{remark}

For integers $1 \leq \ell<n$, Theorem \ref{thm:ThmMain2} and \ref{thm:ThmMain} also apply to the Grassmann variety $X_\ell = \mathrm{Gr}_{\ell,n}(\R)$ of $\ell$-dimensional linear subspaces in the Euclidean space $\R^n$, which represents another novelty. A linear subspace of $\R^n$ is called rational if it admits a basis consisting of elements in $\Q^n$. As in Schmidt's paper \cite{Schmidt67}, we use the Pl\"ucker embedding to define the height $H(v)$ of a rational subspace $v$, and study the approximation of a real subspace chosen at random by rational subspaces. The distance used on $X_{\ell}$ is the usual Riemannian distance and we equip $X_{\ell}$ with the unique probability measure $\sigma_{\ell}$ invariant under rotations. Write $d$ for the dimension $\dim_\R X_{\ell} = \ell(n-\ell)$ and set $\beta_{\ell} = \frac{n}{\ell(n-\ell)}$. By \cite[Th\'eor\`eme~3]{deSaxce22}, for every $c > 0$ and for almost every subspace $x \in X_{\ell}$, the inequality $d(x,v) < c \, H(v)^{-\tau}$ has infinitely (resp. finitely) many solutions $v \in \mathrm{Gr}_{\ell,n}(\Q)$, if $\tau \leq \beta_{\ell}$ (resp. $\tau > \beta_{\ell}$). 

\begin{thm} \label{thm:Grasssmannian}
Fix integers $1 \leq \ell < n$ and let $X_{\ell} = \mathrm{Gr}_{\ell,n}(\R)$ be the Grassmann variety of $\ell$-dimensional subspaces in $\R^n$. Let $\tau \in (0,\beta_{\ell}]$ and $c > 0$. Then there exists an explicit constant $\varkappa > 0$ such that for $\sigma_{\ell}$-almost every subspace $x \in X_{\ell}$, as $T \rightarrow + \infty$,
\begin{equation*} \label{eq:Cor2_ThmMain}
\cN_{c,\tau}(x, T) \sim \varkappa \, c^{\ell(n-\ell)} \, \int_{1}^T y^{(\beta_{\ell}-\tau) \, \ell (n-\ell)} \, \tfrac{d y}{y}.
\end{equation*} 
\end{thm}

\subsection{Method} \label{sec:Method}
We sketch here the method of proof of Theorem \ref{thm:ThmMain2}, which was inspired by the work of Huang and de Saxc\'e~\cite{dSH25} on the local distribution of rational points on flag varieties and uses a celebrated counting method due to Eskin-McMullen~\cite{EM93} and Duke-Rudnick-Sarnak~\cite{DRS93}.

\medskip
Let $\chi$ be a dominant $\Q$-weight of $\bP$ and let $\pi_{\chi} : \bG \rightarrow \GL(\bV_{\chi})$ be the corresponding representation which is strongly rational over $\Q$ (see Section \ref{sec:Height} for the definition and details). In particular, there exists a unique $\bP$-invariant line $\mathbf{D}_{\chi}$ defined over $\Q$ in $\bV_{\chi}$ on which $\bP$ acts via the character $\chi$: 
\[
\forall \, p \in \bP, \forall \, \bm{v} \in \mathbf{D}_{\chi}, \quad \pi_{\chi}(p) \bm{v} = \chi(p) \bm{v}.
\]
Fix a non-zero vector $\bm{e}_\chi \in \mathbf{D}_{\chi}(\Q)$. Thus, we can realize the set of real points $X = \bX(\R)$ as a subvariety of $\mathbb{P}(V_\chi)$, where $V_\chi = \bV_{\chi}(\R)$. Fixing a rational lattice $\bV_\chi(\Z)$ of $V_\chi$ and a $K$-invariant norm $\|\cdot \|$ on $V_\chi$ allows us to define a height function on $\mathbb{P}(\bV_\chi)(\Q)$; by restriction, we then obtain a height function on $\bX(\Q)$. To explain the method, we make the simplifying assumption that $-\bm{e}_{\chi}$ belongs to the $G$-orbit of $\bm{e}_{\chi}$. We call $\cV_{\chi} = G\, \bm{e}_\chi$ the (punctured) cone over $X$; this is the set of all non-zero vectors $ \bm{v} \in V_\chi$ such that the corresponding line $[\bm{v}]$ is in $X$. 
Denote by $\cL_{\chi}$ the set of primitive elements in $\bV_\chi(\Z)$ that are contained in $\cV_{\chi}$.

\medskip
We then translate the problem of counting rational approximations to $x$ of bounded height in $X$ to counting primitive lattice points within a growing family $(k_x \cE_{\psi}(T))_{T \geq 1}$ in the cone $\cV_{\chi}$, where $k_x \in K$ is any element such that $x = k_x P$. In fact, we find that
\[
\cN_\psi(x,T) = \frac{1}{2} \, \# (\cL_{\chi} \cap k_x \cE_{\psi}(T)). 
\]
Hence we are reduced to provide an asymptotic formula for $\# (k_x^{-1}\cL_{\chi} \cap \cE_{\psi}(T))$ for almost every $x \in X$. The stabilizer $\G$ in $G$ of $\bV_\chi(\Z)$ is an arithmetic subgroup of $G$. Let us write $\cE_T$ for $\cE_{\psi}(T)$.

\medskip
\begin{figure}[htbp]
\includegraphics[scale=2.5]{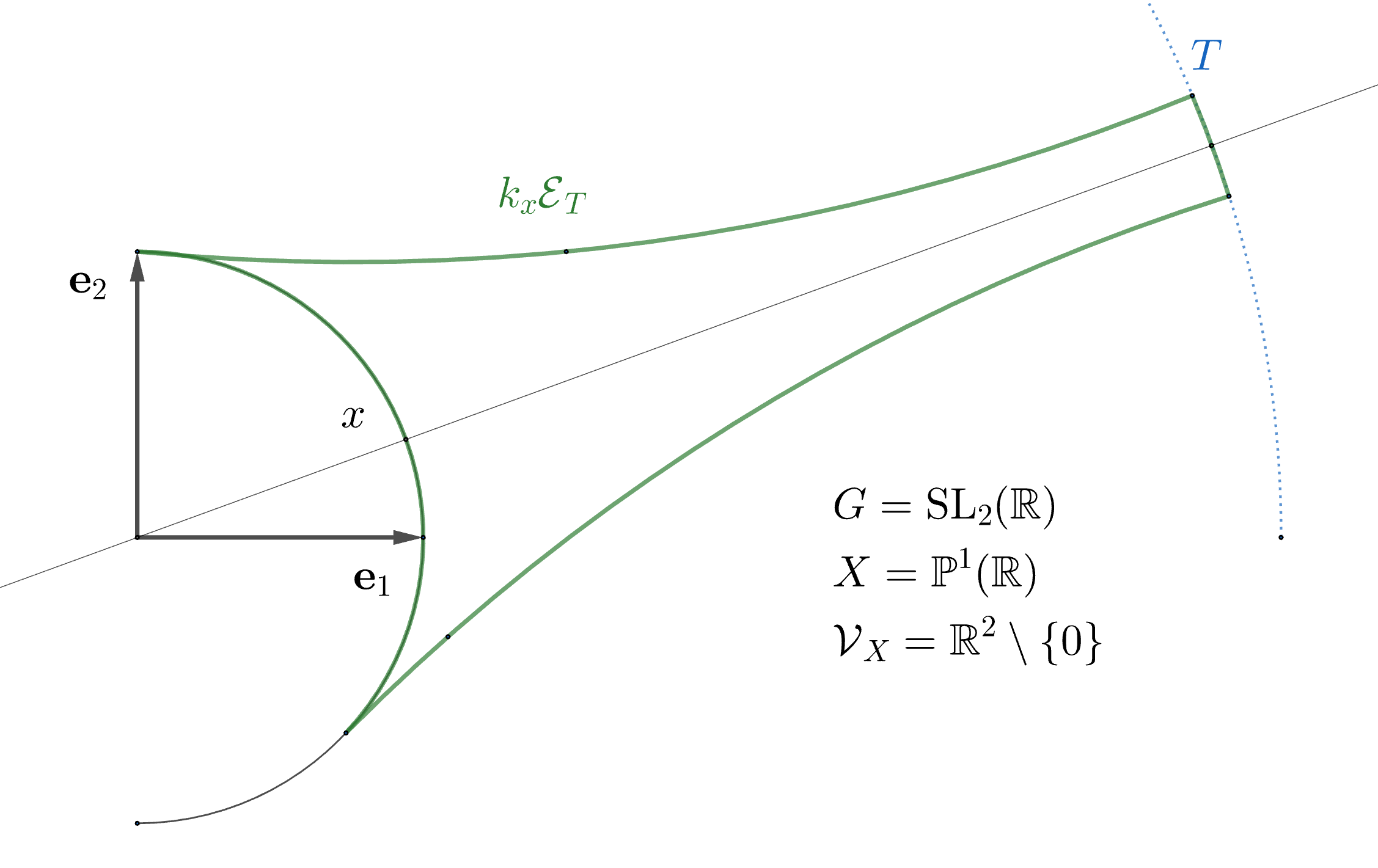}
\caption{The set $k_x \cE_T$ for the group $G = \SL_2(\R)$, $X = \mathbb{P}^1(\R)$, $\cV_{\chi} = \R^2 \setminus \{0\}$, and the set $\cL_{\chi} = \mathcal{P}(\Z^2)$ of primitive elements in $\Z^2$.}
\label{fig:3}
\end{figure}
  
We first note that the set $\cE_T$ is typically not well-rounded (see Definition \ref{def:Well-Rounded}); this term was introduced in \cite{EM93}, it refers to the regularity property of a set to be almost invariant under the action of a small ball centered at the identity in $G$, and it allows for asymptotic lattice point counts with an error term. The idea is then to decompose $\cE_T$ dyadically  
\[
\cE_T = \bigsqcup_{j \geq 0} \cF_{T_j} \quad \text{ where, for all $j\geq 0$, } \cF_{T_j} = \cE_{T /2^{j}} \setminus \cE_{T /2^{j+1}},
\]
and to apply to each $\cF_{T_j}$ an element $a_{y_j}$ of a certain diagonal subgroup $A$ of $G$ in order to obtain a well-rounded set $\cB_{T_j} = a_{y_j} \cF_{T_j}$. 

\medskip
\begin{figure}[htbp]
\includegraphics[scale=0.7]{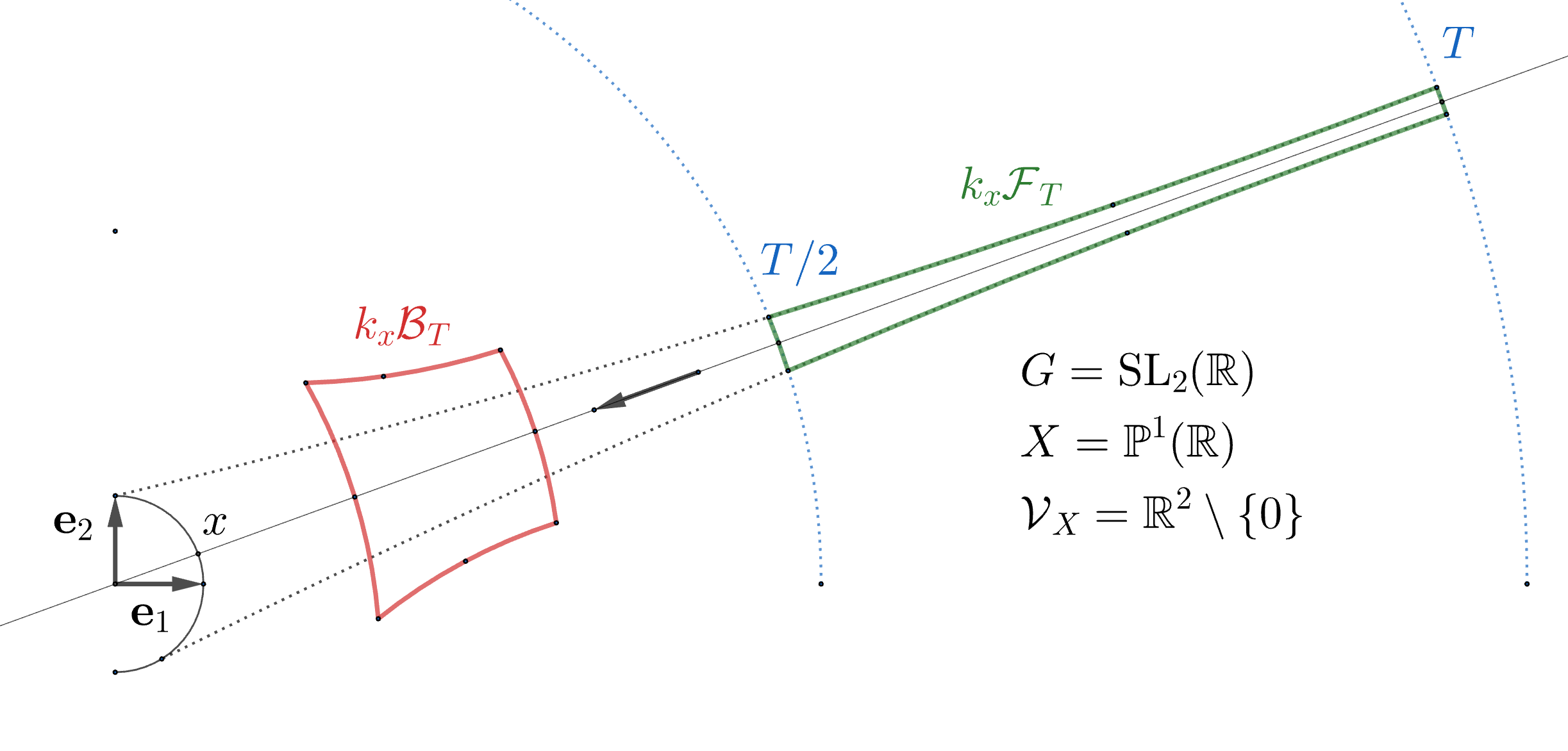}
\caption{The action of $a_{y} = \begin{pmatrix} y^{-1/2} & \\ & y^{1/2}
\end{pmatrix}$ with $y > 1$ on $\cV_{\chi} = \R^2 \setminus \{0\}$ contracts the line through $\bm{e}_\chi = \bm{e}_1$ and expands the line through $\bm{e}_2$. We will show that for a suitable time $y_T > 0$, the set $\cB_{T} = a_{y_T} \cF_{T}$ is well-rounded and the lattice $a_{y_T} k_x^{-1} \Z^2$ is not too distorted for $\sigma_X$-almost every $x \in X$. This suggests that $\# \left (a_{y_T} k_x^{-1} \cL_{\chi} \cap \cB_T \right )$ is approximately given by the volume of $\cB_T$.}
\label{fig:4}
\end{figure}

This leads us to the expression
\[
\# (k_x^{-1}\cL_{\chi} \cap \cE_T) = \sum_{j \geq 0} \# (a_{y_j} k_x^{-1}\cL_{\chi} \cap \cB_{T_j}).
\]
Now, using the exponential mixing property of the diagonal subgroup $A$ and the fact that $\cL_{\chi}$ is a finite union of $\G$-orbits, we show that there exist exists an explicit constant $\varkappa_1 > 0$ and $\varepsilon > 0$ such that for all $T \geq 1$ large and all $g \in G$ whose norm is bounded by $T^{\varepsilon}$, one can estimate $\# \left (g \cL_{\chi} \cap \cB_T \right )$ in terms of the volume of $\cB_T$: 
\[
\# \left (g \cL_{\chi} \cap \cB_T \right ) = \varkappa_1 \, m_{\cV_{\chi}}(\cB_T) \left ( 1 + O( T^{-\varepsilon}) \right ).
\]
On the other hand, using the ergodicity of the action of $A$ on the space of lattices and reduction theory, we can show that the first minimum of $a_{y_j} k_x^{-1}\bV_{\chi}(\Z)$ is typically not affected too much. Putting everything together, we can thus estimate $\# (k_x^{-1} \cL_{\chi} \cap \cE_T)$ in terms of the volume of the set $\cE_T$. Providing the necessary volume estimates then concludes the proof of Theorem \ref{thm:ThmMain2}. 

\subsection{Notation and conventions}
For two positive quantities $A$ and $B$, we use the notation $A\lesssim B$  {or $A=O(B)$} to mean that there is a constant $C > 0$ such that $A\leq CB$, and we use subscripts to indicate the dependence of the constant on parameters. If $A\lesssim B$ and $B\lesssim A$, then we write $A\asymp B$. 

\vspace{5mm}
\textbf{Acknowledgments}. 
I would like to thank Nicolas de Saxc\'e for proposing this problem, for numerous patient discussions about the method of proof and for feedback concerning earlier versions of this paper. I am also grateful to Wooyeon Kim and Gabriel Pallier for sharing with me crucial insights that contributed to the proofs of Lemmas \ref{lem:Equi2} and \ref{lem:Riemann}, respectively. Finally, I am grateful to the anonymous referee for carefully reviewing this paper and providing many valuable suggestions.

\section{Notation and preliminary results} \label{sec:prel}
In this section, we recall the general setting and assumptions from the introduction, equip the set of rational points $\mathbf{X}(\Q)$ with a height function $H_\chi$, specify the distance used on the flag variety $X$, and record some preliminary results from representation theory and lattice reduction.

\subsection{Flag variety and standing assumptions} \label{sec:Flag-Variety}
Let $\bX$ be a generalized flag variety defined over $\Q$, obtained as the quotient $\bX = \bG / \bP$ of a connected semisimple algebraic $\Q$-group $\bG$ by a parabolic $\Q$-subgroup $\bP$. We recall that we assume that the unipotent radical of $\bP$ is abelian. In particular, the $\Q$-rank of $\bX$ is $1$, or equivalently, $\bP$ is a maximal parabolic $\Q$-subgroup of $\bG$. Note that the universal cover $\widetilde{\bG}$ of $\bG$ is a connected simply connected semisimple $\Q$-group and the covering map $\widetilde{\pi} : \widetilde{\bG} \rightarrow \bG$ is a $\Q$-isogeny. Then $\widetilde{\bP} = \widetilde{\pi}^{-1}(\bP)$ is a maximal parabolic $\Q$-subgroup of $\widetilde{\bG}$ and the implied morphism $\widetilde{\bG} / \widetilde{\bP} \rightarrow \bG / \bP$ is an isomorphism defined over $\Q$. Thus, without loss of generality, we can assume that $\bG$ is simply connected.

\medskip
We let $\bP_0$ be a minimal parabolic $\Q$-subgroup contained in $\bP$ and we let $\bT$ be a maximal $\Q$-split $\Q$-torus of $\bG$ contained in $\bP_0$. Let $\Phi$, $\Delta$ and $\{\lambda_{\alpha}\}_{\alpha\in \Delta}$ be the set of roots of $\bG$ relative to $\bT$, the set of simple roots for the ordering associated to $\bP_0$, and the relative fundamental $\Q$-weights (see \cite[Section 3]{MSG14}), respectively. By the maximality of $\bP$, there exists a unique simple root $\alpha \in \Delta$, that we fix throughout this paper, such that $\bP = \bP_{\Delta \setminus \{\alpha\}}$ is the standard parabolic $\Q$-subgroup of $\bG$ associated to the subset $\Delta \setminus \{\alpha\}$ of $\Delta$ (see \cite[Section 11.7]{Borel69}). 

\medskip
We can moreover assume, without loss of generality, that $\bG$ is almost $\Q$-simple. In fact, write $\bG = \bG_1 \cdots \bG_n$ as an almost direct product of connected almost $\Q$-simple $\Q$-groups $\bG_j$, for $1 \leq j \leq n$, and assume that the root $\alpha$ belongs to the set of roots of $\bG_1$ relative to the maximal $\Q$-split $\Q$-torus $(\bG_1 \cap \bT)^\circ$. Then, for all $2 \leq j \leq n$, we have $\bG_j \subset \bP$ and we may assume that $\bG$ is almost $\Q$-simple.

\subsection{Strongly rational representations} \label{sec:Reps}
We briefly recall a few facts from representation theory and refer the reader to \cite[Section 12]{BT65} for details.

\medskip
Let $\bB$ be a Borel subgroup of $\bG$ and suppose that $\bS$ is a maximal $\Q$-torus of $\bG$ contained in $\bB$. Let $\pi : \bG \rightarrow \GL(\bV)$ be an irreducible representation. There exists a unique line $\bD$ in $\bV$ that is invariant under the Borel subgroup $\bB$. In particular, $\bS$ acts on this line through a character $\widetilde{\chi} \in X^*(\bS)$ called the \emph{dominant} or \emph{highest weight} of $\pi$: $\forall \, s \in \bS, \, \forall \, \bm{v} \in \bD$, we have $\pi(s) \bm{v} = \widetilde{\chi}(s) \bm{v}$. The orbit under $\bG$ of the line $\bD$ is called the \emph{cone of the representation $\pi$} and we denote it by $\cC_{\pi}$. The stabilizers of the lines in $\cC_{\pi}$ form a class of parabolic subgroups of $\bG$, denoted by $\cP_{\pi}$, that are conjugated to each other. The representation $(\pi,\bV)$ is called \emph{strongly rational over $\Q$} if it is defined over $\Q$ and if $\cP_{\pi}$ contains a parabolic subgroup $\bP$ which is defined over $\Q$.

\medskip
We suppose now that the maximal $\Q$-split $\Q$-torus $\bT$ is contained in $\bS$ and that the minimal parabolic $\Q$-subgroup $\bP_0$ contains $\bB$. Let $\pi : \bG \rightarrow \GL(\bV)$ be a representation which is strongly rational over $\Q$. Then the stabilizer $\bP$ of the unique $\bB$-invariant line $\bD$ in $\bV$ is a parabolic $\Q$-subgroup containing $\bP_0$. The \emph{weights of $\pi$} are the characters $\mu \in X^*(\bS)$ for which there exists a non-zero $\bm{v} \in \bV$ such that, $\forall \, s \in \bS$, we have $\pi(s) \bm{v} = \mu(s) \bm{v}$. The \emph{$\Q$-weights of $\pi$} are the restrictions to $\bT$ of the weights of $\pi$; the \emph{dominant} or \emph{highest $\Q$-weight of $\pi$}, denoted by $\chi$, is the restriction of $\widetilde{\chi}$ to $\bT$. This character determines $\pi$ up to rational equivalence. Hence we denote from now on $\pi$ by $\pi_{\chi}$ and $\bV$ by $\bV_{\chi}$. By \cite[Proposition~12.13]{BT65}, the highest $\Q$-weight $\chi$ of a representation $\pi_{\chi} : \bG \rightarrow \GL(\bV_{\chi})$ which is strongly rational over $\Q$ is a linear combination $\chi = \sum_{\alpha \in \Delta} c_{\alpha} \, \lambda_{\alpha}$ with $c_{\alpha} \in \N$ of the relative fundamental weights $\{\lambda_{\alpha}\}_{\alpha \in \Delta}$ and, since $\bG$ is assumed to be simply connected, every linear combination $\sum_{\alpha \in \Delta} c_{\alpha} \, \lambda_{\alpha}$ with $c_{\alpha} \in \N$ is the highest $\Q$-weight of a representation $\pi : \bG \rightarrow \GL(\bV)$ which is strongly rational over $\Q$.

\medskip
For every $\Q$-weight $\mu \in X^*(\bT)$ of a representation $\pi_{\chi} : \bG \rightarrow \GL(\bV_{\chi})$ which is strongly rational over $\Q$, let 
\[
\bV^{\mu} = \{\bm{v} \in \bV_{\chi} : \forall \, t \in \bT, \, \pi_{\chi}(t) \bm{v} = \mu(t) \bm{v} \}.
\]
This is a $\Q$-subspace of $\bV_{\chi}$. It is known that $\bV^{\chi}$ is one-dimensional, that $\bV_{\chi}$ is the direct sum of the linear subspaces $\bV^{\mu}$:
\begin{equation} \label{eq:Direct-Sum}
    \bV_{\chi} = \bigoplus_{\mu} \bV^{\mu},
\end{equation}
and that every $\Q$-weight of $\pi_{\chi}$ has the form
\begin{equation} \label{eq:Q-weights}
\mu = \chi - \sum_{\alpha \in \Delta} c_{\alpha}(\mu) \alpha \quad \text{ with $c_{\alpha}(\mu) \in \N$.}
\end{equation}

\subsection{Height on $\bX(\Q)$} \label{sec:Height}
Let $\chi \in X^*(\bP)_{\Q}$ be a non-zero element such that its restriction to $\bT$ is a dominant $\Q$-weight of $\bG$ and let $\pi_{\chi} : \bG \rightarrow \GL(\bV_\chi)$ be the associated representation which is strongly rational over $\Q$ (see Section \ref{sec:Reps}). In particular, there exists a unique $\bP$-invariant line $\bD_{\chi}$ in $\bV_{\chi}$, which is defined over $\Q$ and on which $\bT$ acts through $\chi$. Moreover, since the parabolic $\Q$-subgroup $\bP = \bP_{\Delta \smallsetminus \{\alpha\}}$ is maximal, there exists an integer $n_{\alpha} \in \N^\times$ such that $\chi = n_{\alpha} \lambda_{\alpha}$. For all $g \in \bG$ and $\bm{v} \in \bV_{\chi}$, we abbreviate $\pi_\chi(g) \bm{v}$ as $g \bm{v}$. Fix a non-zero vector $\bm{e}_\chi \in \mathbf{D}_{\chi}(\Q)$ and write $x_0 = [\bm{e}_\chi]$ for the corresponding point in the projective space $\mathbb{P}(\bV_{\chi})$. Then, for every $p \in \bP$, we have $p \bm{e}_\chi = \chi(p) \bm{e}_\chi$ and
\[
\bP = \{g \in \bG : g x_0 = x_0\}.
\]
We identify the orbit $\bG \, x_0$ with $\bX = \bG/\bP$ via the map $\iota_\chi : g \bP \mapsto g x_0$. Let $K$ be a maximal compact subgroup of $G = \bG(\R)$. Let $T = \bT(\R)^\circ$ be the identity component of $\bT(\R)$ in the real topology. Fix a rational basis $(\bm{v}_i)_{i \in I}$ of $V_\chi = \bV_{\chi}(\R)$ and a Euclidean inner product $\langle \cdot, \cdot \rangle$ on $V_{\chi}$ for which the action of $K$ is unitary, that of $T$ is self-adjoint, and such that the basis $(\bm{v}_i)_{i \in I}$ is orthonormal. Moreover, by Equation \eqref{eq:Direct-Sum}, we can assume that this basis consists of eigenvectors for the action of $T$. We denote the implied norm on $V_\chi$ by $\| \cdot \|$. We assume that $\bm{v}_1 = \bm{e}_\chi$. This gives us a height function $H$ on $\mathbb{P}(\bV_\chi)(\Q)$ by $H([\bm{v}])=\|\bm{v}\|$, where $\bm{v}$ is a primitive vector in the lattice $\bV_{\chi}(\Z) = \oplus_{i \in I} \Z \bm{v_i}$ representing $[\bm{v}]$. This height function is well defined, since a primitive vector in $\bV_{\chi}(\Z)$ representing $[\bm{v}]$ is uniquely determined up to multiplication by $\pm 1$. Using the embedding $\iota_\chi$, one then obtains a height function $H_\chi$ on $\bX(\Q)$, which is given by,
\[
\forall \, v \in \bX(\Q), \quad H_\chi(v) = H(\iota_\chi(v)).
\]

\subsection{Measure on $X$}
We equip the maximal compact subgroup $K$ with the Haar probability measure $\mu_K$. By the Iwasawa decomposition $G = KP,$ the group $K$ acts transitively by left multiplication on the set of real points $X = \bX(\R)$. We equip $X$ with the pushforward $\sigma_X$ of $\mu_K$ along the orbital map $k \mapsto k x_0$; this is the unique $K$-invariant probability measure on $X$.

\subsection{Distance on $X$} \label{sec:CC} 
Let us equip the space $X$ with a Riemannian distance that is compatible with the Euclidean structure on $V_\chi$. Let $\mathbb{S} = \{\bm{x} \in V_\chi : \| \bm{x}\| = 1 \}$ be the unit sphere in $V_\chi$, viewed as a Riemannian submanifold of $V_\chi$. The $K$-equivariant projection map $\mathbb{S} \rightarrow \mathbb{P}(V_\chi)$, $\bm{v} \mapsto [\bm{v}]$, being a smooth local diffeomorphism, induces a $K$-invariant Riemannian metric on $\mathbb{P}(V_\chi)$, and by restriction also on $X$. We observe that the associated Riemannian measure  equals $\vol_{\mathrm{R}}(X) \, \sigma_X$, where $\vol_{\mathrm{R}}(X)$ denotes the total Riemannian volume of $X$.
We denote the induced Riemannian distance on $X$ by $d(\cdot, \cdot)$. Let us denote by $\kg$ the Lie algebra of $G$. For every $\beta \in \Phi$, we let $\kg_{\beta} = \{Y \in \kg : \forall \, t \in T, \, \Ad(t) Y = \beta(t) Y \}$ and  $\mathfrak{z} = \{Y \in \kg : \forall \, t \in T, \, \Ad(t) Y = Y \}$. Then, we have the root space decomposition
\[
\kg = \mathfrak{z} \oplus \left ( \bigoplus_{\beta \in \Phi} \kg_{\beta} \right ). 
\]
Let $\Phi^+$ be the set of positive roots and, as before, let $\alpha \in \Delta$ be the unique root such that $\bP = \bP_{\Delta \setminus \{\alpha\}}$ is the standard parabolic $\Q$-subgroup associated to the subset $\Delta \setminus \{\alpha\}$. Let $\langle \Delta \setminus \{\alpha\} \rangle^-$ be the set of all roots that can be expressed as a negative integer linear combination of elements in $\Delta \setminus \{\alpha\}$. Then the Lie algebras of $P$ and of the unipotent subgroup $U^-$ opposite to $P$ are, respectively,
\[
\mathfrak{p} = \mathfrak{z} \oplus \left ( \bigoplus_{\beta \in \Phi^+ \cup \langle \Delta \setminus \{\alpha\} \rangle^-} \kg_{\beta} \right ) \quad \text{and} \quad \ku^- = \bigoplus_{\beta \in \Phi, \, - \alpha \prec \beta} \kg_{\beta},
\]
where $- \alpha \prec \beta$ means that $- \alpha$ occurs in the support of $\beta$. Observe that $\mathfrak{g} = \ku^- \oplus \mathfrak{p}$. Let $\phi : G \rightarrow X$ be the projection map $\phi(g) = g x_0$ and let $D_1 \phi : \kg \rightarrow T_{x_0} X$ be its derivative at the identity $1 \in G$. Observe that $\mathrm{ker} D_1 \phi = \mathfrak{p}$. In particular, if we denote the restriction of this projection map from $G$ to $U^-$ still by $\phi$, then $D_1 \phi : \ku^- \rightarrow T_{x_0} X$ is a linear isomorphism. We equip $\ku^-$ with a Euclidean structure for which this isomorphism is an isometry and we denote the implied norm on $\ku^-$ by $\|\cdot \|_{\ku^{-}}$. It will be convenient for us to relate the distance on $X$ to the one on the Lie algebra $\ku^-$.

\begin{lemma} \label{lem:Riemann} 
For all $u \in \ku^-$,
\begin{equation} \label{eq:Distance_Estimate}
d(x_0, \exp(u) x_0) = \|u\|_{\ku^{-}} + O\bigl (\|u\|_{\ku^{-}}^2 \bigr ),
\end{equation}
where $\exp : \ku^- \rightarrow U^-$ is the exponential map. 
\end{lemma}

\begin{proof}
It suffices to prove \eqref{eq:Distance_Estimate} for all $u \in \ku^-$ sufficiently close to the origin. Let $\| \cdot \|_{x_0}$ the implied norm on $T_{x_0} X$. For all $Y \in T_{x_0}X$ close to the origin, the Riemannian exponential map at $x_0$, that we denote $\exp_{x_0} : T_{x_0} X \rightarrow X$, satisfies: $d(x_0, \exp_{x_0}(Y)) = \| Y \|_{x_0}$. Using the triangle inequality and that the derivative $D_1 \phi$ is an isometry yields, for all $\, u \in \ku^-$ close to the origin,
\begin{align*} \label{eq:Distance_Estimate}
d(x_0, \exp(u) x_0) 
&= \|u\|_{\ku^{-}} + O \left ( d(\exp_{x_0}(D_1 \phi(u)), \exp(u) x_0 )\right ).
\end{align*} 
Moreover, since the Riemannian exponential map $\exp_{x_0}$ is locally bi-Lipschitz, it suffices to show that the map $\Phi = \exp_{x_0}^{-1} \circ \phi \circ \exp$, defined on a neighborhood of the origin in $\ku^-$, satisfies:
\[
\Phi(u) = D_1 \phi(u) + O\bigl (\|u\|_{\ku^{-}}^2 \bigr ).
\]
Since $\Phi(0) = 0$, by Taylor's theorem,
\[
\Phi(u) = D_0 \Phi(u) + O\bigl (\|u\|_{\ku^{-}}^2 \bigr ).
\]
By the chain rule, $D_0 \Phi = D_{x_0} \exp_{x_0}^{-1} \circ D_1 \phi \circ D_0 \exp = D_1 \phi$, as required.
\end{proof}

\subsection{Diophantine exponent} \label{sec:Diophantine}

As in \cite[D\'efinition~2.4.1]{deSaxce20}, we define the \emph{Diophantine exponent $\beta_{\chi}(x)$ relative to $\chi$} of a point $x \in X = \bX(\R)$ by 
\[
\beta_{\chi}(x) = \inf \{\beta > 0 : \exists \, c > 0, \, \forall \, v \in \bX(\Q), \, d(x,v) > c \, H_{\chi}(v)^{-\beta}\}.
\]
By \cite[Theorem~2.4.5]{deSaxce20}, this Diophantine exponent is constant almost everywhere on $X$: there exists a positive rational number $\beta_{\chi} = \beta_{\chi}(X) > 0$ such that for $\sigma_X$-almost every $x \in X$, we have $\beta_{\chi}(x) = \beta_{\chi}.$ We refer to $\beta_{\chi}$ as the \emph{Diophantine exponent of $X$ relative to $\chi$}. By \cite[Th\'eor\`emes~2.4.5 et 3.2.1]{deSaxce20}, for every $c > 0$ and for almost every $x \in X$, the inequality
\begin{equation} \label{eq:DiophantineExponent}
d(x,v) < c \, H_\chi (v)^{-\tau}
\end{equation}
has infinitely (resp. finitely) many solutions $v \in \bX(\Q)$, if $\tau \leq \beta_\chi$ (resp. $\tau > \beta_\chi$).

\subsection{Cone over $X$, coordinates and measures} \label{sec:LightCone}
By abuse of notation, we shall refer to the orbit $\cV_{\chi} = G \, \bm{e}_\chi$ in $V_{\chi}$, which is contained in the punctured affine cone over $X$, as the \emph{cone over $X$}. Let $\cL_{\chi}$ be the set of primitive elements in $\bV_\chi(\Z)$ that are contained in $\cV_{\chi}$. Let $\G$ be the stabilizer in $G$ of $\bV_\chi(\Z)$; this is an arithmetic subgroup of $G$ that stabilizes $\cL_{\chi}$. We normalize the Haar measure $\mu_{G}$ on $G$ so that the induced $G$-invariant measure on the quotient $\Omega$, which we denote by $\mu_{\Omega}$, is a probability measure. 

\smallskip
We briefly recall a few facts from the structure theory of parabolic subgroups and refer the reader to \cite[Section 11.7]{Borel69} for details. By the maximality of $\bP$, there exists a unique simple root $\alpha \in \Delta$ such that $\bP = \bP_{\Delta \setminus \{\alpha\}}$ is the standard parabolic $\Q$-subgroup associated with the subset $\Delta \setminus \{\alpha\}$ of $\Delta$. Define the one-dimensional $\Q$-split $\Q$-torus $\bA = \left ( \bigcap_{\beta \in \Delta \setminus \{\alpha\}} \ker \beta \right )^\circ$ and let $\cZ(\bA)$ be its centralizer in $\bG$. The parabolic $\bP$ is the semidirect product $\bP = \cZ(\bA) \ltimes \bU$, where $\bU$ is its unipotent radical. Let $\bM$ be the identity component of the intersection of the kernels of the $\Q$-characters of $\cZ(\bA)$. By \cite[Proposition~10.7, (b)]{Borel69}, we have $X^*(\bM)_{\Q} = \{1\}$ and $\cZ(\bA)$ is an almost direct product of $\bM$ and $\bA$. Hence $\bP$ has a Langlands decomposition 
\begin{equation} \label{eq:ParabolicDecomposition}
\bP = \bM \, \bA \, \bU.
\end{equation} 
Denote by $\bL$ the stabilizer in $\bG$ of the vector $\bm{e}_\chi$ and note that its identity component (with respect to the Zariski topology) is the semi-direct product $\bL^{\circ} = \bM \ltimes \bU$ and hence does not admit any non-trivial $\Q$-characters. Thus, by \cite[Theorem 9.4]{BHC62}, the Lie group $L = \bL(\R)$ is unimodular and  $\G_L = \G \cap L$ is a lattice in $L$. Let $\mu_L$ be the Haar measure on $L$ normalized so that the induced $L$-invariant measure $\mu_{L/\G_L}$ is a probability measure. We identify the cone $\cV_{\chi}$ with the homogeneous space $G / L$ via the orbital map $g L \mapsto g \bm{e}_{\chi}$.  In particular, by \cite[Theorem 2.51]{F15}, there exists a unique $G$-invariant Radon measure $m_{\cV_{\chi}} = \mu_{G/L}$ on $\cV_{\chi} = G /L$ such that
\begin{align} \label{eq:Folland-Unfolding}
\forall \, f \in C_c(G), \quad \int_G f(g) \, \dd \mu_G(g) &= \int_{G /L} \int_{L} f(g l) \, \dd \mu_{L}(l) \, \dd \mu_{G/L}(g L).
\end{align}

\smallskip
By abuse of notation, we let $A$ be the identity component (with respect to the real topology) of the group of real points $\bA(\R)$. Write $\kt$ for the Lie algebra of $\bT(\R)$ and $Y_\alpha \in \kt$ for the unique element satisfying
\begin{equation} \label{eq:DiagonalElementDef}
\alpha(Y_\alpha) = -1 \quad \text{and} \quad \beta(Y_\alpha) = 0 \quad \text{for all } \beta \in \Delta \setminus \{\alpha\}.
\end{equation}
Note that we have 
\begin{equation} \label{eq:Parametrization}
A = \left \{ a_y = \exp(\log(y) Y_\alpha) : y > 0 \right \}.
\end{equation} 
We recall that $\beta_{\chi}$ denotes the Diophantine exponent of $X$ relative to $\chi$ (see Section \ref{sec:Diophantine}). By \cite[Th\'eor\`eme~2.4.5]{deSaxce20}, we have 
\begin{equation} \label{eq:Relation-Chi-Exponent}
\chi(Y_\alpha) = -\frac{1}{\beta_\chi}
\end{equation}
Hence the action of $a_y$ on $\bm{e}_\chi$ is given by
\begin{equation} \label{eq:a_y-Action-Exponent}
a_y \bm{e}_\chi = y^{-\tfrac{1}{\beta_\chi}} \bm{e}_\chi.
\end{equation}

Let $\rho : \kt \rightarrow \R_+^\times$ be the sum of all the positive roots with multiplicities counted, that is, $\rho = \sum_{\beta \in \Phi^+} \beta,$ where $\Phi^+$ is the set of positive roots of $\bG$ relative to $\bT$ with respect to $\bP_0$. Let $d = \dim X$ be the dimension of $X$ as a real manifold. Then, by the definition of $Y_{\alpha}$ and the proof of \cite[Th\'eor\`eme 2.4.7]{deSaxce20}, we have
\begin{equation} \label{eq:Value-of-Rho}
\rho(Y_{\alpha}) = \sum_{\beta \in \Phi^+} \beta(Y_{\alpha}) = - d.
\end{equation}

Next, we claim that $G$ can be decomposed as $G = K \, A \, L$. In fact, by \cite[Section 11.19, Equation (1)]{Borel69}, we have $G = K \, P = K \, P^\circ$. Moreover, since $P^\circ \subseteq A \, L \subseteq P$, the claim follows. Let us now decompose the Haar measure on $G$ according to the decomposition $G = K \, A \, L$. The modular function of $A \, L$ is the restriction to $A \, L$ of the modular function $\Delta_P$ of $P$. Since $L$ is a closed normal unimodular subgroup of $P$, we have $\Delta_P |_{L} = 1$. By \cite[Section 4, p. 540]{Knapp23}, for every $a \in A$, we have $\Delta_P(a) = \rho(a)$. Moreover, by Equation \eqref{eq:Value-of-Rho} and using the parametrization \eqref{eq:Parametrization} of $A$, for every $y > 0$, we have $\Delta_P(a_y) = \rho(a_y) = y^{-d}$. Let $\mu_K$ be the Haar probability measure on $K$. On $A$ we consider the Haar measure given by the push-forward of the measure $\dd y / y$ on $\R_+^\times$ via the map $y \mapsto a_y$. Therefore, by applying \cite[Theorem 8.32]{Knapp23}, first to the decomposition $G = K \, (A \, L)$, then to the decomposition $A \, L$, there exists a normalizing constant $\omega_0 > 0$ such that the Haar measure $\mu_G$ of $G$ is given by 
\begin{align}\label{equ:Haar1}
\dd \mu_{G}(g)= \omega_0 \, y^{-(d+1)}  \, \dd \mu_{L}(l) \, \dd y \, \dd \mu_K(k).
\end{align}
Moreover, let $K_L = K \cap L$ and let $\sigma$ be the pushforward of the measure $\mu_K$ on $K$ to $K_L$ via the map $k \mapsto k \, K_L$. Using the decomposition $G = K \, A \, L$, the map of $(K/K_L) \times A$ to $\cV_{\chi} = G / L$ given by $(k \, K_L, a_y) \mapsto k a_y \bm{e}_{\chi}$ is a homeomorphism. Let us now show that the measure $m_{\cV_{\chi}}$ is given by
\begin{align}\label{eq:Measure-On_Tilde-X}
\dd m_{\cV_{\chi}}(k a_y \bm{e}_\chi) = \omega_0 \, y^{-(d+1)} \, \dd y \, \dd \sigma(k).
\end{align}
By the uniqueness (up to scalars) of the $G$-invariant Radon measure on $\cV_{\chi} = G/L$ (see, for instance, \cite[Theorem 2.51]{F15}), it suffices to show the $G$-invariance of the measure on the right hand side in \eqref{eq:Measure-On_Tilde-X}. So fix $g_0 \in G$ and $f \in C_c(\cV_{\chi})$. By \cite[Proposition 2.50]{F15}, there exists $\phi \in C_c(G)$ such that, for every $g \in G$,
\[
f(g \bm{e}_{\chi}) = \int_{L} \phi(g l) \, \dd \mu_L(l).
\]
Hence, using the fact that $\rho(a_y) = y^{-d}$, that $\sigma$ is the pushforward from $K$ to $K/K_L$ of $\mu_K$, the measure description in \eqref{equ:Haar1} of $\mu_G$ and its invariance, we have
\begin{align*}
\int_{\cV_{\chi}} &f(g_0 k a_y \bm{e}_{\chi}) \, \omega_0 \, y^{-(d+1)} \, \dd y \, \dd \sigma(k)  = \int_{\cV_{\chi}} \int_{L} \phi(g_0 k a_y l) \, \dd \mu_L(l) \, \omega_0 \, y^{-(d+1)} \, \dd y \, \dd \sigma(k)  \\
&= \int_{K/K_L} \int_{A} \int_{L} \phi(g_0 k a_y l) \, \dd \mu_L(l) \, \omega_0 \, y^{-(d+1)} \, \dd y \, \dd \sigma(k)  \\
&= \int_{K} \int_{A} \int_{L} \phi(g_0 k a_y l) \, \omega_0 \, y^{-(d+1)} \, \dd \mu_{L}(l) \, \frac{\dd y}{y} \, \dd \mu_K(k) \\
&= \int_{G} \phi(g_0 g) \, \dd \mu_G(g) = \int_{G} \phi(g) \, \dd \mu_G(g).
\end{align*}
By the same argument, we have 
\[
\int_{\cV_{\chi}} f(k a_y \bm{e}_{\chi}) \, \omega_0 \, y^{-(d+1)} \, \dd y \, \dd \sigma(k) = \int_{G} \phi(g) \, \dd \mu_G(g).
\]
This shows that
\[
\int_{\cV_{\chi}} f(g_0 k a_y \bm{e}_{\chi}) \, \omega_0 \, y^{-(d+1)} \, \dd y \, \dd \sigma(k) = \int_{\cV_{\chi}} f(k a_y \bm{e}_{\chi}) \, \omega_0 \, y^{-(d+1)} \, \dd y \, \dd \sigma(k)
\]
and, since $g_0 \in G$ was arbitrary, the measure on the right-hand side in \eqref{eq:Measure-On_Tilde-X} is $G$-invariant, as required.

\subsection{Sobolev norms} \label{sec:Sobolev}
We recall that $\Omega$ denotes the homogeneous space $G/\G$. In this section, following \cite[Sections 3.7]{EMV09}, for every $\ell \in \N^\times$, we define a degree $\ell$ Sobolev norm $\cS_{\ell}$ on the space $C_c^\infty(\Omega)$ of compactly supported smooth functions on $\Omega$. Let us first introduce some notation and record relevant facts from \cite[Sections 3.6, 3.7]{EMV09}. We fix an Euclidean norm $\|\cdot \|_{\kg}$ on the Lie algebra $\kg$ of $G$ such that $\|[\bm{u}, \bm{v}] \|_{\kg} \leq \|\bm{u} \|_{\kg} \| \bm{v} \|_{\kg}$ and choose a rational $\Ad(\G)$-stable lattice $\kg_\G \subset \kg$ satisfying $[\kg_\G, \kg_\G] \subseteq \kg_\G$. 
For $x \in \Omega = G / \G$, we then set
\[
\mathrm{ht}(x) = \sup \, \left \{ \|\Ad(g) \bm{v}\|_{\kg}^{-1} : x = g\G, \, \bm{v} \in \kg_\G \smallsetminus \{\bm{0}\} \right \}.
\]
Let $d_G(\cdot, \cdot)$ be a right-invariant Riemannian metric on $G$ and, for every $r > 0$, denote by $B_G(r)$ the ball in $G$ with center $1 \in G$ and radius $r$. For every $x \in \Omega$, denote by $\mathrm{inj}(x)$ the injectivity radius at $x$. By lattice reduction theory, there exist constants $c_3 > 0$ and $\kappa_0 > 0$ such that, for every $x \in \Omega$, the map $B_G(c_3 \mathrm{ht}(x)^{-\kappa_0}) \rightarrow \Omega$ given by $g \mapsto g x$ is injective (see \cite[Section 3.6, Equation (3.7)]{EMV09}).

\medskip
For every lattice $\Delta$ in $\kg$, we denote by $\lambda_1(\Delta)$ its first minimum, that is, the length $\|\bm{v}\|$ of the shortest non-zero vector $\bm{v}$ in $\Delta$. In particular, we have that $\mathrm{ht}(g\G) = \lambda_1(\Ad(g)\kg_\G)^{-1}$. Now, choose an orthonormal basis for $\kg$ and, for every $\ell \in \N^\times$, define  the $L^2$-Sobolev norm $\cS_{\ell}$ on $C_c^\infty(\Omega)$ by
\[
\cS_{\ell}(f)^2 = \sum_{\deg \cD \leq \ell} \|\mathrm{ht}(x)^\ell \, \cD f \|_{L^2(\Omega)}^2,
\]
where the sum is taken over all $\cD \in U(\kg)$, the universal enveloping algebra of $\kg$, which are monomials in the chosen basis of degree $\leq \ell$.

\begin{lemma} \label{lem:Bump-Function-Cusp}
For every $\ell \in \N^\times$ there exists $D \in \N^\times$ such that for every $x \in \Omega$ and $\delta \in (0, \mathrm{inj}(x))$, there exists a non-negative function $\phi \in C_c^\infty(\Omega)$ satisfying
\begin{equation} \label{eq:Bump-Function-Cusp}
\mu_{\Omega}(\phi) = 1, \quad \supp(\phi) \subseteq B_{\Omega}(x, \delta), \quad \text{and} \quad \cS_{\ell} (\phi) \lesssim \delta^{-D} \mathrm{ht}(x)^{\ell}.
\end{equation}
\end{lemma}
\begin{proof}
Let $\ell \in \N^\times$, $x \in \Omega$ and $\delta \in (0, \mathrm{inj}(x))$. By \cite[Lemma~2.4.7]{KM96}, there exists a non-negative function $\varphi \in C_c^\infty(G)$, only depending on $\ell$ and $r$, satisfying
\begin{equation} \label{eq:Bump-Function-Cusp-G}
\mu_{G}(\varphi) = 1, \,\,\, \supp(\varphi) \subseteq B_{G}(1, \delta/10), \,\,\, \text{and} \,\,\,  \Big (\sum_{\deg \cD \leq \ell} \| \cD \varphi \|_{L^2(G)}^2 \Big )^{1/2} \lesssim \delta^{-(\ell + \dim G /2)}.
\end{equation}
We define, by convolution, $\phi = \varphi \star \big ( c^{-1} \, \mathbbm{1}_{B_{\Omega}(x,\delta/10)} \big)$ with $c = \mu_{\Omega}(B_{\Omega}(x,\delta/10))$. Then, we have $\supp(\phi) \subseteq  B_{\Omega}(x,\delta)$ and, using Fubini's theorem, $\mu_{\Omega}(\phi) = 1$. Finally, using the Cauchy-Schwarz inequality and the fact that, for every differential operator $\cD \in U(\kg)$, we have $\cD \phi =  (\cD \varphi) \star \big ( c^{-1} \, \mathbbm{1}_{B_{\Omega}(x,\delta/10)} \big)$ , yields
\begin{align*}
\cS_{\ell}(\phi)^2 &\lesssim \mathrm{ht}(x)^{2 \ell} \, \sum_{\deg \cD \leq \ell} \| \cD \phi \|_{L^2(\Omega)}^2 \\ &= \mathrm{ht}(x)^{2 \ell} \, c^{-2} \, \sum_{\deg \cD \leq \ell} \int_{\Omega} \Big | \int_G (\cD \phi)(g^{-1}) \mathbbm{1}_{B_{\Omega}(x,\delta/10)} (g y) \, \dd \mu_G(g) \Big |^2 \dd \mu_{\Omega}(y) \\
&\lesssim \mathrm{ht}(x)^{2 \ell} \, c^{-2} \, \sum_{\deg \cD \leq \ell} \int_{\Omega} \int_{G} \big | \cD \phi \big |^2 \, \dd \mu_G \, \int_{B_G(\delta)} \mathbbm{1}_{B_{\Omega}(x,\delta/10)} (g y) \, \dd \mu_G(g) \, \dd \mu_{\Omega}(y) \\
&= \mathrm{ht}(x)^{2 \ell} \, c^{-2} \, \Big ( \sum_{\deg \cD \leq \ell} \| \cD f \|_{L^2(G)}^2 \Big ) \int_{B_G(\delta)} \, \int_{\Omega} \mathbbm{1}_{B_{\Omega}(x,\delta/10)} (g y) \, \dd \mu_{\Omega}(y) \, \dd \mu_G(g) \\
&\lesssim \mathrm{ht}(x)^{2 \ell} \, \delta^{-2 (\ell + \dim G /2)}.
\end{align*}
This completes the proof of Lemma \ref{lem:Bump-Function-Cusp} with $D = \ell + \dim G /2$.
\end{proof}

\subsection{Lattice reduction with respect to $\bM$} \label{sec:Lattice}
Let $\bL$, $\bM$, $\bA$ and $\bU$ be as in Section \ref{sec:LightCone}. In the proof of Theorem \ref{thm:ThmMain2}, we will need the effective equidistribution of translates of certain $L$-orbits. We note that $\bL^\circ$ is the semi-direct product of the connected reductive $\Q$-group $\bM$ and the unipotent radical $\bU$ of $\bP$. We will deduce the required equidistribution result (see Lemma \ref{lem:Equi2}) by using the effective equidistribution of expanding translates of (horospherical) $U$-orbits (see \cite[Proposition 2.4.8.]{KM96}) together with lattice reduction with respect to $M$. The latter will be the content of this section.

\medskip
By \cite[Theorem 9.4]{BHC62}, since $\bM$ does not admit any non-trivial $\Q$-characters, we have that $\G_M = \Gamma \cap M$ is a lattice in $M$. We note that $K_M = K \cap M$ is a maximal compact subgroup of $M$, $\bP_0 \cap \bM$ is a minimal parabolic $\Q$-subgroup of $\bM$ and $\bS = (\bT \cap \bM)^\circ$ is a maximal $\Q$-split $\Q$-torus of $\bM$ (see \cite[Section 11.7]{Borel69}). Then $[\Delta \smallsetminus \{\alpha\}]$, the set of roots that are linear combinations of elements of $\Delta \smallsetminus \{\alpha\}$, is the root system of $\bM$ relative to $\bS$. Let $\bU_0$ denote the unipotent radical of $\bP_0 \cap \bM$ and let $\bM_0$ be the maximal connected $\Q$-anisotropic $\Q$-subgroup of the centralizer $\cZ(\bS)$ of $\bS$ in $\bM$. We let $\ks$ denote the Lie algebra of $S = \bS(\R)$. For every $t > 0$, let 
\[
\ks_t = \{X \in \ks : \forall \beta \in \Delta \smallsetminus \{\alpha\}, \, \beta(X) \leq t \}
\]
The negative Weyl chamber of $\ks$ is given by
\[
\ks^- = \{X \in \ks : \forall \beta \in \Delta \smallsetminus \{\alpha\}, \, \beta(X) < 0 \}.
\]
We set $S^- = \exp (\ks^-)$ and note that there exists a compact neighborhood $Q$ of the identity in $S$ such that $\exp(\ks_t) \subset Q \, S^-$. Denote by $M_0$ and $U_0$ the groups of real points of $\bM_0$ and $\bU_0$, respectively. By \cite[Theorem 13.1]{Borel69}, there exists $t > 0$, a compact neighborhood $\bm{\omega}$ of the identity in $M_0 U_0$ and a finite set $C \subset \bM(\Q)$, such that the Siegel set $\mathfrak{S} = \mathfrak{S}_{\tau, \bm{\omega}} = K_M \, \exp(\ks_t) \, \bm{\omega}$ satisfies $M = \mathfrak{S} \, C \, \G_M$. For every $\delta \in (0,1)$, we define
\begin{equation} \label{eq:Siegel-Set-M}
\mathfrak{S}(\delta)  = \{m \in \mathfrak{S} : \lambda_1(\Ad(m) \kg_\G) < \delta \},
\end{equation}
and we denote by $\mathfrak{S}(\delta)^c$ the complementary subset of $\mathfrak{S}(\delta)$ in $\mathfrak{S}$. Let $\mu_M$ be the Haar measure on $M$. There exists $c_1 > 0$ such that, for every $\delta \in (0,1)$, we have 
\begin{equation} \label{eq:Siegel-Cusp-Bound}
\mu_{M}(\mathfrak{S}(\delta)) \,\lesssim \, \delta^{c_1}.
\end{equation}
In fact, by \cite[Lemma 12.2]{Borel69}, the union $\cup_{s  \in \exp(\ks_t)} s \bm{\omega} s^{-1}$ is relatively compact, and, for all $k \in K_M \, Q$, $s \in S^-$, $n \in \bm{\omega}$, $c \in C$, and $\gamma \in \G_M$, we have 
\begin{equation} \label{eq:Comparable1}
\lambda_1(\Ad(k s n c \gamma) \kg_\G) \, \asymp \, \lambda_1(\Ad(s) \kg_\G).
\end{equation}
Moreover, since $S$ is contained in the torus $T$ and the roots $\Phi$ are the non-trivial $\Q$-weights in the adjoint representation of $G$, 
\begin{equation} \label{eq:Comparable2}
\forall \, s \in S, \quad \lambda_1(\Ad(s) \kg_\G) \asymp \min_{\beta \in \Phi} \, \beta(s).   
\end{equation}
Choose a character $\chi' \in X^*(\bS)$ such that 
\[
\forall \, s \in S^-, \quad \chi'(s) \leq \min_{\beta \in \Phi} \, \beta(s).
\]
Let $\rho' : \ks \rightarrow \R_+^\times$ be the sum of all the positive $S$-roots with multiplicities counted. By \cite[Theorem 8.32]{Knapp23}, we can decompose the measure $\mu_M$ on $M$ according to the decomposition $M = K_M \, M_0 \, S \, U_0$ and deduce that, for some constant $C' \geq 1$, 
\begin{align*}
\mu_{M}(\mathfrak{S}(\delta)) \, &\lesssim \, \int_{S^-} \mathbbm{1}_{\{\chi'(s) \, \leq \, C' \, \delta \}} \rho'(s) \, \dd \mu_S(s) \\
&= \, \int_{\ks^-} \mathbbm{1}_{\{\chi'(s) \, \leq \, \ln (C' \, \delta) \}} \, e^{\rho'(s)} \, \dd s.  
\end{align*}
Note that the $\Q$-rank of $\bM$ is $r = \mathrm{rank}_{\Q}(\bM) = \mathrm{rank}_{\Q}(\bG) - 1$. Let $(Y_i)_{1 \leq i \leq r}$ be the dual basis of the set of simple roots of $\bM$ relative to $\bS$. This last integral is the integral of an exponential function over a convex polytope in $\ks$, and is therefore comparable to the maximum of the function on this set, multiplied by a factor corresponding to the dimension of the face on which this maximum is attained. This maximum is attained at one of the vertices of the convex set, which are the points $A_i = \frac{\ln (C' \, \delta)}{\chi'(Y_i)} Y_i$ with $i = 1, \dots, r$. Hence, for all $\delta \in (0,1)$ sufficiently small, we have
\[
\int_{\ks^-} \mathbbm{1}_{\{\chi'(s) \, \leq \, \ln (C' \, \delta) \}} \, \rho'(s) \, \dd s \asymp (C'\delta)^{a} |\ln (C'\delta) |^{b-1} \lesssim \delta^{a/2}
\]
where 
\[
a = \min_{1 \leq i \leq r} \frac{\rho'(Y_i)}{\chi'(Y_i)} \quad \text{and} \quad  b = \# \{ 1 \leq i \leq r : \frac{\rho'(Y_i)}{\chi'(Y_i)} = a \}.
\]
This proves the claim in \eqref{eq:Siegel-Cusp-Bound} with $c_1  = a/2$.

Moreover, for every $m \in \mathfrak{S}(\delta)^c$, we have 
\begin{equation} \label{eq:Siegel-Adjoint-Bound}
\|\Ad(m) \| \lesssim \delta^{-1}.
\end{equation}
Indeed, by \eqref{eq:Comparable1} and \eqref{eq:Comparable2}, for all $m = k s n c \gamma$ with $k \in K_M \, Q$, $s \in S^-$, $n \in \bm{\omega}$, $c \in C$, and $\gamma \in \G_M$, we have 
\begin{align*}
\|\Ad(m) \| &\asymp \|\Ad(s) \| \asymp \max_{\beta \in \Phi} \beta(s) = \left ( \min_{\beta \in \Phi} \beta(s) \right )^{-1} \\
&\asymp \lambda_1(\Ad(s) \kg_\G)^{-1} \asymp \lambda_1(\Ad(m) \kg_\G)^{-1} \leq \delta^{-1}.
\end{align*}

\subsection{Exponential mixing}
Finally, we will need the following theorem, which follows from the results of \cite[Section 3.4]{KM99} and \cite[Section 4.1-4.3]{EMMV20} applied to $\bG$, $\bT$ and $\G$ as in Sections \ref{sec:Flag-Variety} and \ref{sec:LightCone}.
\begin{thm} \label{thm:Exponential-Mixing}
Let $Y$ be an element in the Lie algebra of $\bT(\R)$ such that, for every projection $\pi_i : \mathfrak{g} \to \mathfrak{g}_i$ onto a $\mathbb{Q}$-simple factor, $\pi_i(Y) \neq 0$. For every $y \in \R_+^\times$, let $a_y' = \exp(\ln(y) \, Y)$. Then there exist constants $\ell \in \mathbb{N}$ and $C', c' > 0$ such that for all $\phi_1, \phi_2 \in C_c^\infty(\Omega)$ and all $0 < y \leq 1$,
\[
\left | \int_\Omega \phi_1(a_y' x) \, \phi_2(x) \, \dd \mu_\Omega(x)
\, - \, 
\mu_\Omega(\phi_1) \, \mu_\Omega(\phi_2) \right |
 \, \leq \, 
C'\, y^{c'} \, \|\Upsilon^\ell \phi_1\|_2\, \|\Upsilon^\ell \phi_2\|_2,
\]
where $\Upsilon$ denotes the differential operator
\[
\Upsilon = 1 - \sum_i Y_i^2,
\]
where $(Y_i)$ is an orthonormal basis of the Lie algebra $\mathfrak{k}$ of the maximal compact subgroup $K$ of $G$.
\end{thm}

\section{Diophantine approximation and counting lattice points} \label{sec:Reduction}
In this section, we translate the problem of counting rational approximations of bounded height in $X$ to the problem of counting primitive lattice points in a certain family of growing sets in the Euclidean space $V_\chi$. 

\subsection{Counting lattice points}
Let $\psi : \R_{+}^{\times} \rightarrow \R_{+}^{\times}$ be a decreasing function. For every $x \in X$ and $T \geq 1$, we define the counting function 
\[
\cN_{\psi}(x, T) = \# \, \big \{v \in \bX(\Q) : d(x,v) < \psi(H_{\chi}(v)), \, H_{\chi}(v) < T \big \}.
\]
Moreover, for every $T \geq 1$, we associate to $\psi$ the set
\begin{equation} \label{eq:Definition_E}
\cE_{\psi}(T) = \big \{ \bm{v} \in \cV_{\chi} : d(x_0, [\bm{v}]) < \psi(\|\bm{v}\|), 1 \leq \|\bm{v}\| < T \big \}.
\end{equation}
Fix a section $\ks : X  \rightarrow K$ of the orbital map $K \to X$ sending $k$ to $kx_0$. Given $x \in X$, we shall write $k_x = \ks(x)$.
As the following lemma shows, estimating the counting function $\cN_\psi(x,T)$ amounts to counting lattice points in the increasing family $\{\cE_T\}_{T \geq 1}.$ Let $[K \cap P : K \cap L] \in \{1,2\}$ be the index of $K \cap L$ in $K \cap P$.
\begin{lemma} \label{lem:Reduction}
For every $x \in X$ and $T \geq 1$,
\begin{equation} \label{eq:Reduction}
\cN_\psi(x,T) = [K \cap P : K \cap L]^{-1} \, \# \left (k_x^{-1} \cL_{\chi}  \cap \cE_{\psi}(T) \right ).
\end{equation}
\end{lemma}

\begin{proof}
It suffices to show that $\cN_\psi(x,T) = \# \left (\cL_{\chi}  \cap k_x \cE_{\psi}(T) \right )$. We first note that
\[
k_x \cE_{\psi}(T) = \left \{ \bm{v} \in \cV_{\chi} : d(x, [\bm{v}]) < \psi(\|\bm{v}\|), 1 \leq \|\bm{v}\| < T \right \}.
\]
Now a rational point $v = g [\bm{e}_\chi] \in \bX(\Q)$ satisfies $d(x, v) < \psi(H_\chi(v))$ and $1 \leq H_\chi(v) < T$ if and only if any of the primitive vectors $\bm{v} \in \cL_{\chi}$ representing $v$ satisfies $d(x, [\bm{v}]) < \psi(\|\bm{v}\|)$ and $1 \leq \|\bm{v}\| < T$. This finishes the proof of the lemma.
\end{proof}

Our main Theorem \ref{thm:ThmMain2} follows from the following statement, which provides an asymptotic formula in terms of the measure of the set $\cE_{\psi}(T)$ and whose proof will be given in Sections~\ref{sec:WellRoundedness} and \ref{sec:Mixing}. Recall that the measure $m_{\cV_{\chi}}$ on $\cV_{\chi}$ has been constructed in \eqref{eq:Measure-On_Tilde-X} and that $d = \dim X$ denotes the dimension of $X$ as a real manifold.

\begin{proposition} \label{prop:Counting}
Let $\psi : \R_{+}^{\times} \rightarrow \R_{+}^{\times}$ be a non-increasing function. Suppose that there exist $\tau \in (0,\beta_\chi)$ and $C > 1$ such that $C^{-1} y^{-\tau} \leq \psi(y) \leq C y^{-\tau}$ for all sufficiently large $y \geq 1$. Then there exist $\varepsilon_1 > 0$ and an explicit constant $\varkappa_1 > 0$ such that for $\sigma_X$-almost every $x \in X$, as $T \rightarrow + \infty$, we have
\[
\# \left (k_x^{-1} \cL_{\chi}\cap \cE_{\psi}(T) \right ) =  \varkappa_1 \, m_{\cV_{\chi}}(\cE_{\psi}(T)) \left ( 1 + O_x ( m_{\cV_{\chi}}(\cE_{\psi}(T))^{-\varepsilon_1})\right).
\]
\end{proposition}

\begin{proof} [Proof of Theorem \ref{thm:ThmMain2}]
Let $\psi : \R_{+}^{\times} \rightarrow \R_{+}^{\times}$ be as in Theorem \ref{thm:ThmMain2}. By Lemma \ref{lem:Reduction} and Proposition \ref{prop:Counting}, it only remains to provide the desired measure estimate for $\cE_{\psi}(T)$. Let $\sigma$ be the pushforward of the Haar probability measure $\mu_K$ on $K$ to the quotient $K / K_L$ along the quotient map $k \mapsto k K_L$. We recall from Section \ref{sec:LightCone} that the map of $(K/K_L) \times A$ to $\cV_{\chi} = G / L$ given by $(k \, K_L, a_y) \mapsto k a_y \bm{e}_{\chi}$ is a homeomorphism and, in these coordinates, 
the measure $m_{\cV_{\chi}}$ is given by
\[
\dd m_{\cV_{\chi}}(k a_y \bm{e}_\chi) = \omega_0 \, y^{-(d+1)} \, \dd \sigma(k) \, \dd y.
\]
By Equation \eqref{eq:a_y-Action-Exponent}, the action of $a_y$ on $\bm{e}_\chi$ is given by
\[
a_y \bm{e}_\chi = y^{-\tfrac{1}{\beta_\chi}} \bm{e}_\chi.
\]
Using the fact that the norm $\| \cdot \|$ on $V_{\chi}$ is $K$-invariant and $\|\bm{e}_{\chi} \| = 1$, the set $\cE_{\psi}(T)$ can now be described as
\[
\cE_{\psi}(T) = \left \{ k a_y \bm{e}_\chi \in \cV_{\chi} : d(x_0,k x_0)  < \psi(y^{-\frac{1}{\beta_\chi}}) , \, 1 \leq y^{-\frac{1}{\beta_\chi}} < T \right\}.
\]
We denote by $B_X(r)$ the ball of radius $r > 0$ centered at $x_0 \in X$. We compute
\begin{align} \label{eq:Error}
m_{\cV_{\chi}}(\cE_{\psi}(T)) &=   \int_{T^{-\beta_\chi}}^{1} \int_{K / K \cap L} \mathbbm{1}_{B_X( \psi(y^{-\frac{1}{\beta_\chi}})))}(k x_0) \, \omega_0 \, \dd \sigma(k) \frac{d y}{y^{d+1}} \notag \\
&= \int_{T^{-\beta_\chi}}^{1} \sigma_X(B_X( \psi(y^{-\frac{1}{\beta_\chi}}) )) \, \omega_0 \, \frac{d y}{y^{d+1}}.
\end{align} 
Let $\vol_{\mathrm{R}}(X)$ be the total Riemannian volume of $X$. As is true for any Riemannian manifold \cite[Theorem~3.1]{G74}, for every $r > 0$, we have
\begin{equation} \label{eq:RiemannianVolume}
\sigma_X(B_X(r)) = \vol_{\mathrm{R}}(X)^{-1} \, \frac{\pi^{d/2}}{\Gamma(\tfrac{d}{2} + 1)} \, r^d + O(r^{d+2}),  
\end{equation}
where $\frac{\pi^{d/2}}{\Gamma(\frac{d}{2} + 1)}$ is the volume of the unit ball in the Euclidean space $\R^d$. Using the change of variable $y \mapsto y^{-\frac{1}{\beta_\chi}}$, the assumption $\psi(y) \asymp y^{-\tau}$, and defining, for every $T \geq 1$, the integral $\Psi(T) = \int_{1}^T \psi(y)^{d} y^{\beta_\chi d} \frac{d y}{y}$, we have
\begin{align} \label{align:Volume-Computation}
m_{\cV_{\chi}}(\cE_{\psi}(T)) &= \beta_\chi \,  \int_{1}^{T} \sigma_X(B_X(\psi(y))) y^{\beta_\chi d} \, \omega_0 \, \frac{d y}{y} \\ \notag
&= \beta_\chi \, \frac{\pi^{d/2}}{\Gamma(\tfrac{d}{2} + 1)} \, \vol_{\mathrm{R}}(X)^{-1} \, \omega_0 \int_{1}^{T} \psi(y)^d y^{\beta_\chi d} \frac{d y}{y} + O \left (\int_{1}^{T} y^{(\beta_\chi - \tau) d - 2\tau} \frac{d y}{y} \right )\\\notag
&= \beta_\chi \, \frac{\pi^{d/2}}{\Gamma(\tfrac{d}{2} + 1)} \, \vol_{\mathrm{R}}(X)^{-1} \, \omega_0 \, \Psi(T) + O \left ( T^{(\beta_\chi - \tau) d - 2\tau} \right )\\ \notag
&= \beta_\chi \, \frac{\pi^{d/2}}{\Gamma(\tfrac{d}{2} + 1)} \, \vol_{\mathrm{R}}(X)^{-1} \, \omega_0 \, \Psi(T) \left ( 1 + O(\Psi(T)^{- \frac{2\tau}{(\beta_\chi - \tau) d}}) \right ),
\end{align} 
Let $\varepsilon_1 > 0$ and $\varkappa_1 > 0$ be as in Proposition \ref{prop:Counting}. This establishes the proof of Theorem \ref{thm:ThmMain2} with 
\begin{equation} \label{eq:Varkappa-Constant}
\varkappa = [K \cap P:K \cap L]^{-1} \, \varkappa_1 \, \beta_\chi \,
\frac{\pi^{d/2}}{\Gamma(\tfrac{d}{2} + 1)}  \, \vol_{\mathrm{R}}(X)^{-1}  \, \omega_0
\end{equation}
and $\varepsilon = \min \{\varepsilon_1, \frac{2\tau}{(\beta_\chi - \tau) d}\}$.
\end{proof}

\section{Diagonal action and well roundedness} \label{sec:WellRoundedness}
To establish the lattice point counting in Proposition \ref{prop:Counting}, we use the exponential mixing property of the flow $A$ and the well-roundedness property of certain sets, in the spirit of the counting method developed by Eskin-McMullen~\cite{EM93} and Duke-Rudnick-Sarnak~\cite{DRS93}. Fix a decreasing function $\psi : \R_{+}^{\times} \rightarrow \R_{+}^{\times}$ as in Theorem \ref{thm:ThmMain2}. However, the associated set $\cE_{\psi}(T)$ is typically not well-rounded with respect to the action of $G$ (see Figure \ref{fig:3}). For every $T \geq 1$, we define 
\[
\cF_{T} = \left \{ \bm{v} \in \cV_{\chi}: d(x_0, [\bm{v}]) < \psi(\|\bm{v}\|), \, \max \{ 1 ,{T}/2 \} \leq \|\bm{v}\| < {T} \right \}.
\]
The idea is to first decompose $\cE_{\psi}(T)$ dyadically:
\[
\cE_{\psi}(T) = \bigsqcup_{j \geq 0} \cF_{T_j} \quad \text{ where } \,  T_j = T /2^{j}, \, j \geq 0.
\]
It will turn out that, for every $T \geq 1$, if we set $y_T = T^{ \tau}$, then the set 
\[
\cB_T = a_{y_T} \cF_T
\]
is well-rounded as defined below. The purpose of this section is to establish the well-roundedness of the family $\cB = (\cB_T)_{T \geq 1}$.

\medskip
We shall work with the following definition. We recall that $d_{G}(\cdot,\cdot)$ denotes a right-invariant Riemannian metric on $G$ and $B_{G}(r)$ the metric open ball with radius $r > 0$ and center $1 \in G$. Our definition is inspired by that of \cite[Section 3.2]{KY23}, which in turn is inspired by the work of \cite{EM93}.

\begin{definition} \label{def:Well-Rounded}
A family $\cB' = (\cB_T')_{T \geq 1}$ of finite-measure Borel subsets of $\cV_{\chi}$ is said to be \emph{well-rounded}, if there exist constants  $C_1 >  0$, $\delta_0 > 0$, and $T_0 \geq 1$ such that for any $\delta \in (0,\delta_0)$ and $T > T_0$, the sets
\[
\overline{\cB}_{T,\delta}' = \bigcup_{g \in B_G(\delta)} g \cB_T' \quad \text{and} \quad \underline{\cB}_{T,\delta}' = \bigcap_{g\in B_G(\delta)} g \cB_T'
\]
satisfy 
\begin{equation} \label{eq:well-rounded}
m_{\cV_{\chi}}\left ( \overline{\cB}_{T,\delta}' \setminus \underline{\cB}_{T,\delta}' \right ) \leq C_1 \delta m_{\cV_{\chi}}(\cB_T').
\end{equation}
\end{definition}
By Section \ref{sec:Reps} (in particular, Equation \eqref{eq:Direct-Sum}), the vector space $V_\chi$ decomposes as a finite direct sum of $\Q$-weight spaces with respect to the action of $T$. In particular, by restricting this action to $A$ and letting $(\chi_j)_{j \in J}$ be the finite family of induced $\Q$-weights of $A$ and $(V_j)_{j\in J}$ the corresponding $\Q$-weight spaces, explicitly given by 
\[
V_j = \{\bm{v} \in V_\chi : \, \forall \, y \in \R_+^\times, \, a_y \bm{v} = \chi_j(a_y) \bm{v} \},
\]
the vector space $V_\chi$ decomposes as a finite direct sum
\begin{equation} \label{eq:Decompose}
V_\chi = \bigoplus_{j \in J} V_j.
\end{equation}
We can assume that the subspaces $V_j$ are mutually orthogonal. 
We let $\pi^+ : V_\chi \rightarrow V_\chi$ be the orthogonal projection onto $\R\bm{e}_\chi$ and we simply write $\bm{v}^+$ for $\pi^+(\bm{v})$.

\medskip
For any $\delta \in (0,1)$ and $y>0$, we define 
\begin{equation} \label{eq:Psi_PM}
\psi_\delta^{-}(y) = (1 + \delta)^{- 1} \psi((1 +\delta) y) \quad \text{and} \quad \psi_\delta^{+}(y) = (1 + \delta) \psi((1 +\delta)^{- 1} y).
\end{equation}

Using these functions $\psi_\delta^{\pm} : \R_+ \rightarrow \R_+$, we define the sets
\[
\underline{\cF}_{T,\delta} = \left\{ \bm{v} \in \cV_{\chi}: d(x_0, [\bm{v}]) < \psi_{\delta}^{-}(\|\bm{v}\|), \, (1 + \delta)\tfrac{T}{2} < \|\bm{v} \| < (1 + \delta)^{-1} T \right\}
\]
and
\[
\overline{\cF}_{T,\delta} = \left\{ \bm{v} \in \cV_{\chi}: d(x_0, [\bm{v}]) < \psi_{\delta}^{+}(\|\bm{v}\|),\, (1 + \delta)^{-1}\tfrac{T}{2} < \|\bm{v} \| < (1 + \delta) T \right\}.
\]
Then we put 
\begin{equation} \label{eq:DefinitionB_T}
\underline{\cB}_{T,\delta} = a_{y_T} \underline{\cF}_{T,\delta} \quad \text{and} \quad \overline{\cB}_{T,\delta} = a_{y_T} \overline{\cF}_{T,\delta}.
\end{equation}

\begin{lemma} \label{lem:WellRoundedCone}
For all $T \geq 1$ large enough and $\bm{w} \in \cB_T$, we have $\|\bm{w}\| \asymp \|\bm{w}^+\|$.
\end{lemma}

\begin{proof}
For all $T \geq 1$ large enough, every $\bm{v} \in \cF_T$ satisfies that $[\bm{v}] \in X$ is close to $x_0$. Since the map $\phi : \ku^- \rightarrow X$ given by $\phi : u^- \mapsto \exp(u^-) x_0$ parametrizes a neighborhood of $x_0 \in X$, there exists a unique $u_{\bm{v}}^- \in \ku^-$ such that $[\bm{v}] = \exp(u_{\bm{v}}^-) x_0$. Moreover, there exists a unique $y_{\bm{v}} \in \R_+^\times$ such that $\bm{v} = \pm \exp(u_{\bm{v}}^-) a_{y_{\bm{v}}} \bm{e}_{\chi}$. Without loss of generality, we may assume that $\bm{v} = \exp(u_{\bm{v}}^-) a_{y_{\bm{v}}} \bm{e}_{\chi}$.  Using the definition of $a_{y} = \exp(\log (y) Y_\alpha)$ and the fact that $\ku^-$ is abelian, the adjoint action of $a_{y}$ on $\ku^- = T_{x_0} X$ is given by $\Ad(a_{y}) u^- = y u^-$ (see \cite[Lemme~2.4.2]{deSaxce20}). Hence the element $\bm{w} = a_{y_T} \bm{v} \in \cB_T$ can be written as
\[
a_{y_T} \bm{v} = \exp(y_T u_{\bm{v}}^-) a_{y_T} a_{y_{\bm{v}}} \bm{e}_{\chi}.
\]
The elements $y_T u_{\bm{v}}^-$ stay in a fixed compact neighborhood $V$ of the origin in $\ku^-$. In fact, by the estimate \eqref{eq:Distance_Estimate}, we have $\|u_{\bm{v}}^- \|_{\ku^-} \asymp d(x_0, [\bm{v}]) \leq \psi(\|\bm{v}\|) \lesssim T^{-\tau}$, and hence 
\[
y_T \, \|u_{\bm{v}}^- \|_{\ku^-} \lesssim 1.
\]
This implies that 
\begin{equation} \label{eq:Asymp-1}
\| \bm{w} \| \asymp \| a_{y_T} a_{y_{\bm{v}}} \bm{e}_{\chi} \| = \chi( a_{y_T} a_{y_{\bm{v}}} ).
\end{equation}
On the other hand, we observe that for every $u^- \in \ku^-$ and $y \in \R_+^\times$, we have 
\[
a_y \exp(u^-) \bm{e}_{\chi} = \exp(y u^-) \chi(a_y) \bm{e}_{\chi}.
\]
Thus the vector $\exp(u^-) \bm{e}_{\chi}$ gets expanded under the action of $a_y$ at the highest possible rate $\chi(a_y)$ as $y$ tends to zero and we must have $\| (\exp(u^-) \bm{e}_{\chi} )^+ \| > 0$. By the compactness of $V$, we thus have $\| (\exp(u^-) \bm{e}_{\chi} )^+ \| \asymp 1$ for all $u^- \in V$. This yields 
\begin{equation} \label{eq:Asymp-2}
\| \bm{w}^+ \| = \chi(a_{y_T} a_{y_{\bm{v}}}) \| (\exp(y_T u_{\bm{v}}^-) \bm{e}_{\chi} )^+ \| \asymp \chi(a_{y_T} a_{y_{\bm{v}}}). 
\end{equation}
Putting \eqref{eq:Asymp-1} and \eqref{eq:Asymp-2} together completes the proof of Lemma \ref{lem:WellRoundedCone}.
\end{proof}

\begin{lemma} \label{lem:WellRounded} 
There exists $c_1 > 0$ such that for all $T \geq 1$ large enough, $\delta \in (0,1)$ and $g \in B_G(c_1\delta)$,
\begin{equation} \label{eq:verify_well_rounded}
\underline{\cB}_{T,\delta} \subseteq g \cB_T \subseteq \overline{\cB}_{T,\delta}.
\end{equation}  
\end{lemma}

\begin{proof}
We prove the right inclusion in \eqref{eq:verify_well_rounded}; the proof of the inclusion $\underline{\cB}_{T,\delta} \subseteq g\cB_T$ is essentially identical and we omit the details. To simplify the notation, for $g \in G$ and $T \geq 1$, we write $g_T = a_{y_T}^{-1} g a_{y_T}$ and for $\bm{v} \in \cV_{\chi}$ we denote $v = [\bm{v}]$ the corresponding projective point. 

\medskip
Fix $\delta \in (0, 1)$. First, using that $\cB_{T} = a_{y_T} \cF_T$ and $\overline{\cB}_{T,\delta} = a_{y_T} \overline{\cF}_{T,\delta}$, we note that the inclusion $g \cB_T \subseteq \overline{\cB}_{T,\delta}$ holds if and only if $g_T \cF_{T,\delta} \subseteq \overline{\cF}_{T,\delta}$. Therefore, using the definitions of $\cF_{T,\delta}$ and $\overline{\cF}_{T,\delta}$, we need to show that there exists a constant $c_1 > 0$ such that for all $T \geq 1$, $\bm{v} \in \cF_T$, and $g \in B_G(c_1\delta)$, we have
\begin{align*}
(1)&  \quad (1 + \delta)^{-1}\frac{T}{2} \leq \|g_T \bm{v}\| < (1 + \delta) T,\\
(2)& \quad d(x_0, g_T v  ) < \psi_{\delta}^+(\|g_T \bm{v} \|).
\end{align*}
Since each $\bm{v} \in \cF_T$ satisfies $T/2 \leq \|\bm{v} \| < T$, in order to prove $(1)$, it suffices to show 
\begin{equation} \label{eq:1_lem:WellRounded}
(1+\delta)^{-1} \|\bm{v} \| < \|g_T \bm{v}\| < (1+\delta) \|\bm{v} \|.
\end{equation}
We first show the right inequality $\|g_T \bm{v} \| < (1+ \delta) \|\bm{v} \|$ in \eqref{eq:1_lem:WellRounded}. Applying the triangle inequality, one has $\|g_T \bm{v}\| \leq \| \bm{v}\| + \|g_T \bm{v} -  \bm{v}\|$. Thus, we further reduce to showing that for all $g \in  B_G(\delta)$,
\begin{equation} \label{eq:2_lem:WellRounded}
\|g_T \bm{v} - \bm{v} \| \lesssim \delta \|\bm{v} \|.
\end{equation}
Then, by decomposing $\bm{v} = \sum_{j \in J} \bm{v}_j$ into weight vectors according to \eqref{eq:Decompose}, for every $j \in J$, we have
\[
\| g_T \bm{v}_j - \bm{v}_j \| \leq y_T^{\frac{1}{\beta_\chi}} \|g a_{y_T} \bm{v}_j - a_{y_T} \bm{v}_j \| 
\lesssim y_T^{\frac{1}{\beta_\chi}} \delta \|a_{y_T} \bm{v}_j \| 
\lesssim \delta \|\bm{v}^+\|
\]
where for the last inequality we used that $\| a_{y_T} \bm{v}_j \| \lesssim \|a_{y_T} \bm{v}^+\|$, which follows from Lemma \ref{lem:WellRoundedCone}. This implies \eqref{eq:2_lem:WellRounded}, which is the right hand side of \eqref{eq:1_lem:WellRounded}. The proof of the other inequality in \eqref{eq:1_lem:WellRounded} is very similar and we omit the details. 
    
\medskip
Let us now show assertion (2). Using \eqref{eq:1_lem:WellRounded} and the fact that $\psi$ is non-increasing,
\[
\psi_{\delta}^+(\|g_T\bm{v}\|) = (1 + \delta) \psi ((1 + \delta)^{-1} \|g_T\bm{v}\|) \geq (1+\delta) \psi(\|\bm{v} \|).
\]
Thus, it is enough to check that 
\begin{equation} \label{eq:check}
d(x_0, g_T v ) < (1+\delta) \psi(\|\bm{v} \|).
\end{equation}
Using the triangle inequality and the fact that $\bm{v} \in \cF_T$ satisfies $d(x_0, v) < \psi(\|\bm{v}\|)$, 
\[
d(x_0, g_T v ) 
\leq \psi(\|\bm{v}\|) + d(v, g_T v ).
\]
Therefore, proving $(2)$ reduces to showing that, for all $g \in  B_G(\delta)$,
\begin{equation} \label{eq:WellRoundedReduce}
d(v, g_T v ) \lesssim \delta \psi(\|\bm{v} \|). 
\end{equation}
Letting $u_{\bm{v}}^- \in \ku^-$ be such that $v = \exp(u_{\bm{v}}^-) x_0$, one then has 
\begin{align*}
d(v, g_T v) = d(\exp(u_{\bm{v}}^-) x_0, g_T \exp(u_{\bm{v}}^-) x_0) \lesssim d(x_0, \exp(-u_{\bm{v}}^-) g_T \exp(u_{\bm{v}}^-) x_0).
\end{align*}
Conjugating $\exp(u_{\bm{v}}^-)$ by $a_{y_T}$ yields
\[
\exp(-u_{\bm{v}}^-) g_T \exp(u_{\bm{v}}^-) = a_{y_T}^{-1} \exp( - \Ad(a_{y_T}) u_{\bm{v}}^-) g \exp( \Ad(a_{y_T}) u_{\bm{v}}^-) a_{y_T}.
\]
Using the decomposition $\kg = \ku^{-} \oplus \kp$, where $\kp$ is the Lie algebra of $P$, we may write $\exp(-\Ad(a_{y_T}) u_{\bm{v}}^-) g \exp( \Ad(a_{y_T}) u_{\bm{v}}^-) = \exp(u') p$ with $u' \in \ku^-$ and $p \in P$. The distance we want to bound is then
\[
d(v,  g_T v ) \lesssim d(x_0, a_{y_T}^{-1} \exp(u') x_0 ) = d(x_0, \exp(\Ad(a_{y_T}^{-1}) u') x_0 ) \lesssim y_T^{-1} \|u' \|_{\ku^{-}}.
\]
Using the definition $y_T = T^{\tau}$ and that $T^{-\tau} \lesssim \psi(\|\bm{v} \|)$ gives $ y_T^{-1} \lesssim \psi(\|\bm{v}\|)$.
On the other hand, using again that $\|\Ad(a_{y_T}) u_{\bm{v}}^- \|_{\ku^{-}} \lesssim 1 $, one has 
\begin{align*}
\| u' \|_{\ku^{-}} &\lesssim d_G(1, \exp(u')p) 
\\ &= d_G(\exp( \Ad(a_{y_T}) u_{\bm{v}}^-), g \exp( \Ad(a_{y_T}) u_{\bm{v}}^-)) \lesssim d_G(1, g) < \delta,
\end{align*}
as desired.
\end{proof}

\begin{proposition} \label{prop:WellRounded}
The family $(\cB_T)_{T \geq 1}$ is well-rounded.
\end{proposition}

\begin{proof}
Let $c_1 > 0$ be as in Lemma \ref{lem:WellRounded}. We show that there exists $C_1 > 0$ such that for any $T \geq 1$ and for any $\delta \in (0, 1)$, the Borel subsets $\overline{\cB}_{T,\delta},\underline{\cB}_{T,\delta}$ as defined in \eqref{eq:DefinitionB_T} satisfy
\begin{equation*}
\underline{\cB}_{T,\delta} \subseteq \bigcap_{g \in B_G(c_1 \delta)} g \cB_T \subseteq \bigcup_{g \in B_G(c_1 \delta)} g \cB_T \subseteq \overline{\cB}_{T,\delta}.
\end{equation*}
and $m_{\cV_{\chi}}(\overline{\cB}_{T,\delta} \setminus\underline{\cB}_{T,\delta}) \leq C_1 \delta m_{\cV_{\chi}}(\cB_T)$. The inclusion relations follow from Lemma \ref{lem:WellRounded}. Thus, using the $G$-invariance of the measure $m_{\cV_{\chi}}$, it suffices to show the measure bound 
\[
m_{\cV_{\chi}}(\overline{\cF}_{T,\delta} \setminus\underline{\cF}_{T,\delta}) \leq C_1 \delta m_{\cV_{\chi}}(\cF_T).
\]
We have $\overline{\cF}_{T,\delta} \setminus\underline{\cF}_{T,\delta} \subseteq R_{1} \cup R_{2} \cup R_{3}$, where
\[
R_1 = \left \{ k a_y \bm{e}_\chi \in \cV_{\chi}: \begin{aligned} \psi^{-}_{\delta} \left (y^{-\frac{1}{\beta_\chi}} \right ) \leq d(x_0,kx_0)  < \psi^+_{\delta} \left (y^{-\frac{1}{\beta_\chi}} \right ), \\ (1+\delta) \tfrac{T}{2} \leq y^{-\frac{1}{\beta_\chi}} < (1+\delta)^{-1} T  
\end{aligned} \right\},
\]
\[
R_2 = \left \{ k a_y \bm{e}_\chi \in \cV_{\chi}: \begin{aligned} d(x_0,k x_0)  < \psi^+_{\delta} \left (y^{-\frac{1}{\beta_\chi}} \right ), \\ (1+\delta)^{-1} \tfrac{T}{2} \leq y^{-\frac{1}{\beta_\chi}} < (1+\delta) \frac{T}{2}  
\end{aligned} \right\},
\]
and 
\[
R_3 = \left \{ k a_y \bm{e}_\chi \in \cV_{\chi}: \begin{aligned} d(x_0,k x_0)  < \psi^+_{\delta} \left (y^{-\frac{1}{\beta_\chi}} \right ), \\ (1+\delta)^{-1} T \leq y^{-\frac{1}{\beta_\chi}} < (1+\delta) T  
\end{aligned} \right\}.
\]
Let us start with giving an upper bound on the measure of $R_1$. 
By the arguments given in \eqref{align:Volume-Computation}, using the change of variable $y \mapsto y^{-\frac{1}{\beta_\chi}}$, we have
\begin{align*}
m_{\cV_{\chi}}(R_1) &\lesssim  \int_{(1+\delta)\tfrac{T}{2}}^{(1+\delta)^{-1} T} \sigma_X \left (B_X( \psi_\delta^+(y)) \smallsetminus B_X( \psi_\delta^-(y))\right ) y^{\beta_\chi d} \, \frac{\dd y}{y}.
\end{align*}
By \cite[Theorem~3.1]{G74}, for all $0 < r_2 < r_1$ small enough, the volume of the difference of two small balls centered at $x_0$ satisfies 
\begin{equation} \label{eq:RiemannSurface}
\sigma_X(B_X(r_1) \setminus B_X(r_2)) \lesssim r_1^d - r_2^d.
\end{equation}
Using this estimate, we have $m_{\cV_{\chi}}(R_1) \lesssim I_1 - I_2$, where
\[
I_1 = \int_{(1+\delta)\frac{T}{2}}^{(1+\delta)^{-1}T}\psi^+_{\delta}(y)^d y^{\beta_\chi d} \frac{\dd y}{y}, \quad I_2 = \int_{(1+\delta)\frac{T}{2}}^{(1+\delta)^{-1}T} \psi^-_{\delta}(y)^d y^{\beta_\chi d} \frac{\dd y}{y}. 
\]
Using the change of variable $y \mapsto (1+\delta)^{-1} y$ in $I_1$ and rearranging, we obtain $I_1 = I_1^{(1)} - I_1^{(2)}$, where
\begin{align*}
I_1^{(1)} &= (1 + \delta)^{d(1+\beta_\chi)} \int_{\frac{T}{2}}^{T} \psi(y)^d y^{\beta_\chi d} \frac{\dd y}{y}, \\ 
I_1^{(2)} &= (1 + \delta)^{d(1+\beta_\chi)} \int_{(1+\delta)^{-2}T}^{T} \psi(y)^d y^{\beta_\chi d} \frac{\dd y}{y}.
\end{align*}
Similarly, by the change of variable $y \mapsto (1+\delta) y$ in $I_2$ and rearranging, we have $I_2 = I_2^{(1)} - I_2^{(2)}$, where
\begin{align*}
I_2^{(1)} &= (1 + \delta)^{-d(1+\beta_\chi)} \int_{\frac{T}{2}}^{T} \psi(y)^d y^{\beta_\chi d} \frac{\dd y}{y}, \\ I_2^{(2)} &= (1 + \delta)^{-d(1+\beta_\chi)}\int_{\frac{T}{2}}^{(1+\delta)^2 \frac{T}{2}} \psi(y)^d y^{\beta_\chi d} \frac{\dd y}{y}.
\end{align*}
Let us show that $I_1^{(1)} - I_2^{(1)}, I_1^{(2)}, I_2^{(2)} \lesssim \delta m_{\cV_{\chi}}(\cF_T)$. We have
\begin{align*}
I_1^{(1)} - I_2^{(1)} &= \left ( (1 + \delta)^{d(1+\beta_\chi)} - (1 + \delta)^{-d(1+\beta_\chi)} \right ) \, \int_{\frac{T}{2}}^{T} \psi(y)^d y^{\beta_\chi d} \frac{\dd y}{y} \\
&\lesssim \delta \, m_{\cV_{\chi}}(\cF_T).
\end{align*}
Next, using that $\psi(y) \asymp y^{-\tau}$ and evaluating the integral, we have
\begin{align*}
I_1^{(2)} &= (1 + \delta)^{d(1+\beta_\chi)} \int_{(1+\delta)^{-2}T}^{T} \psi(y)^d y^{\beta_\chi d} \frac{\dd y}{y} \\
&\lesssim (1 + \delta)^{d(1+\beta_\chi)} \int_{(1+\delta)^{-2}T}^{T} y^{(\beta_\chi-\tau) d - 1} \, \dd y \\
&\lesssim (1 + \delta)^{d(1+\beta_\chi)} \left ( T^{(\beta_\chi-\tau) d} - ((1+\delta)^{-2}T)^{(\beta_\chi-\tau) d}\right ) \\
&= (1 + \delta)^{d(1+\beta_\chi)} \left ( 1 -((1+\delta)^{-2})^{(\beta_\chi-\tau) d} \right ) T^{(\beta_\chi-\tau) d} \\
&\lesssim \delta \, T^{(\beta_\chi-\tau) d} \asymp \delta \, m_{\cV_{\chi}}(\cF_T). 
\end{align*}
By symmetry, we also have $I_2^{(2)} \lesssim \delta \, m_{\cV_{\chi}}(\cF_T)$, as required. The measure bounds $m_{\cV_{\chi}}(R_2) \lesssim \delta \, m_{\cV_{\chi}}(\cF_T)$ and $m_{\cV_{\chi}}(R_3) \lesssim \delta \, m_{\cV_{\chi}}(\cF_T)$ are shown similarly, and we omit the details.
\end{proof}

\section{From exponential mixing to counting} \label{sec:Mixing}
In this section, we use the exponential mixing property of $A$, an ingredient from geometry of numbers, and the well-roundedness of the family $(\cB_T)_{T \geq 1}$ in order to prove Proposition \ref{prop:Counting}.

\subsection{Equidistribution of translated $L$-orbits} 
The following equidistribution result, derived from the effective equidistribution of translated horospherical orbits (see \cite[Proposition~2.4.8]{KM96}), is probably standard. For the sake of completeness, we include a proof. We recall that $\bL$ denotes the stabilizer in $\bG$ of the vector $\bm{e}_{\chi}$ and that $\Omega$ denotes the homogeneous space $G/\G$. For every $g \in \bG(\Q)$, we consider the stabilizer $\G_L^g = g \G g^{-1} \cap L$ in $L$ of the rational point $x = g \G$ in $\Omega$ and we equip the quotient $L / \G_L^g$ with the unique finite $L$-invariant measure $\mu_{L / \G_L^g}$ induced from the Haar measure $\mu_L$ on $L$.

\begin{lemma} \label{lem:Equi2}
Let $Q$ be a compact subset of $G$. There exist $c > 0$ and $\ell \in \N^\times$, depending only on $G$, such that, for every $g \in \bG(\Q)$, $q \in Q$, $0 < y \leq 1$, and $\phi \in C_c^\infty(\Omega)$,
\begin{equation} \label{eq:Equi3}
\left | \int_{L / \G_L^g} \phi(q a_y l g) \, \dd \mu_{L / \G_L^g} (l) - \mu_{L / \G_L^g}(1) \, \mu_{\Omega}(\phi)  \right | \, \lesssim_{g, G, \G, Q} \, y^{c} \, \mathcal{S}_{\ell}(\phi).
\end{equation}  
\end{lemma}

\begin{proof}
We recall from Section \ref{sec:Lattice} that $\bL^\circ$ is the semi-direct product of the reductive $\Q$-subgroup $\bM$ of $\bG$ and the unipotent radical $\bU$ of $\bP$. Without loss of generality, we may assume that the total volume of $L / \G_L^g$ is $1$. Note that the group $(\bL^\circ)(\R)$ is open and closed, normal, and of finite index in $L = \bL(\R)$, so $(\bL^\circ)(\R) / (\G^g \cap (\bL^\circ)(\R))$ is open and closed in $L / \G_L^g$, and the latter is the union of a finite number of translates $h (\bL^\circ)(\R) / (\G^g \cap (\bL^\circ)(\R))$ with $h \in \cZ(A)$ (since $L \subset P = \cZ(A) \ltimes U$ and $U$ is connected). Hence, we may further assume that $L = (\bL^\circ)(\R)$. In particular, we have $L = M \ltimes U$. Note that $\G_U^g = \G^g \cap U$ and $\G_M^g = \G^g \cap M$ are lattices in $U$ and $M$, respectively. Let us denote by $\mu_{U / \G_U^g}$ and $\mu_{M / \G_M^g}$ the unique $U$ and $M$-invariant probability measures on $U / \G_U^g$ and $M / \G_M^g$, respectively. Let $Q'$ be a compact subset of $G$. We recall that $\bG$ is assumed to be simply connected and almost $\Q$-simple (see Section \ref{sec:Flag-Variety}). Hence, by the definition of $Y_{\alpha}$ in \eqref{eq:DiagonalElementDef}, for every projection $p_i : \kg \rightarrow \kg_i$ onto a $\Q$-simple factor $\kg_i$ of $\kg$, we have $p_i(Y_{\alpha}) \neq 0$. Thus, by Theorem \ref{thm:Exponential-Mixing}, the flow $a_y = \exp( \ln(y) \, Y_{\alpha})$ is exponential mixing on $\Omega$. Therefore, by the proof of \cite[Proposition~2.4.8]{KM96}, there exists a constant $c' > 0$ and $\ell \in \N^\times$, depending only on $G$, such that, for every $g \in \bG(\Q)$, $q \in Q'$, $0 < y \leq 1$, and $\phi \in C_c^\infty(G/\G)$, we have
\begin{equation} \label{eq:Equi4}
\left | \int_{U / \G_U^g} \phi(q a_y u g) \, \dd \mu_{U / \G_U^g}(u) - \mu_{\Omega}(\phi) \right | \, \lesssim_{g, G, \G, Q'} \, y^{c'} \, \mathcal{S}_{\ell}(\phi).
\end{equation} 
By \cite[Corollary 6.4]{BHC62}, the semi-direct product $\G_M^g \, \G_U^g$ has finite index in $\G_L^g$, and we may, without loss of generality, assume that $\G_L^g = \G_M^g \, \G_U^g$. Let us now show that, for every $\phi \in C_c^\infty(\Omega)$, 
\begin{equation} \label{eq:semi-direct-integral}
\int_{L / \G_L^g} \phi(q a_y l g) \, \dd \mu_{L / \G_L^g} (l) = \int_{M / \G_M^g} \int_{U / \G_U^g} \phi(q a_y m u g) \, \dd \mu_{U / \G_U^g} (u) d \mu_{M / \G_M^g}(m).
\end{equation}
We first observe that conjugation by $\gamma \in \G_M^g$ defines a map from $U/\G_U^g$ to itself and sends the $U$-invariant probability measure $\mu_{U/ \G_U^g}$ to a $U$-invariant probability measure, which, by uniqueness, must be $\mu_{U/ \G_U^g}$. Thus, for every $\phi' \in C_c(L/\G_L^g)$, the map $m \mapsto \int_{U / \G_U^g} \phi'( m u) \, \dd \mu_{U / \G_U^g}$ on $M$ defines a well-defined function on $M / \G_M^g$. Note that
\[
\Lambda(\phi') = \int_{M / \G_M^g} \int_{U / \G_U^g} \phi'( m u) \, \dd \mu_{U / \G_U^g} (u) \dd \mu_{M / \G_M^g}(m)
\]
defines a positive $L$-invariant linear functional on $C_c(L / \G_L^g)$. Thus, by the Riesz-Markov-Kakutani representation theorem, there exists a unique constant $c_0 > 0$ such that for every $\phi' \in C_c(L / \G_L^g)$, we have
\[
\int_{L / \G_L^g} \phi'(l) \, \dd \mu_{L / \G_L^g} (l) = c_0 \, \int_{M / \G_M^g} \int_{U / \G_U^g} \phi'(h u) \, \dd \mu_{U / \G_U^g} (u) \dd \mu_{M / \G_M^g}(h).
\]
By Lebesgue's dominated convergence theorem, this equality still holds for $\phi' = 1$. Hence, we have $c_0 = 1$. Since the restriction of $\phi(q a_y \cdot)$ to $L\G^g / \G^g = L/\G_L^g$ is a continuous compactly supported function, equation \eqref{eq:semi-direct-integral} follows. 
\smallskip 
As in Section \ref{sec:Lattice}, we let $\mathfrak{S}$ be a Siegel set of $M$ with respect to $\G_M^g$ and we let $C \subset \bM(\Q)$ be the finite subset such that $M = \mathfrak{S} \, C \, \G_M^g$. Therefore, for every $g \in \bG(\Q)$, $q \in Q$, $0 < y \leq 1$ and $\phi \in C_c^\infty(\Omega)$, 
\begin{align*}
\Bigg | \int_{L / \G_L^g} &\phi(q a_y l g) \, \dd \mu_{L / \G_L^g} (l) - \mu_{\Omega}(\phi)  \Bigg | \\
&= \Bigg | \int_{M / \G_M^g} \int_{U / \G_U^g} \phi(q m a_y u g) \, \dd \mu_{U / \G_U^g} (u) \dd \mu_{M / \G_M^g}(m) - \mu_{\Omega}(\phi) \Bigg | \\
&\lesssim \int_{\mathfrak{S}} \sum_{c \in C} \Bigg | \int_{U / \G_U^g} \phi(q m c a_y u g) \, \dd \mu_{U / \G_U^g} (u) - \mu_{\Omega}(\phi) \Bigg | \, \dd \mu_{M}(m).
\end{align*}
We recall from Section \ref{sec:Lattice} that, for every $\delta \in (0,1)$, we defined
\[
\mathfrak{S}(\delta) = \{m \in \mathfrak{S} : \lambda_{1}(\Ad(h) \kg_{\G} ) < \delta\}.
\]
By \eqref{eq:Siegel-Cusp-Bound}, there exists a constant $c_1 > 0$ such that, for all $\delta \in (0,1)$ small enough, we have $\mu_{M}(\mathfrak{S}(\delta)) \lesssim \delta^{c_1}$. In particular, using that $|\phi| \leq \| \phi \|_{L^\infty(\Omega)}$, we have 
\begin{equation} \label{eq:Proof-Equi-Margulis}
\int_{\mathfrak{S}(\delta)} \sum_{c \in C} \Bigg | \int_{U / \G_U^g} \phi(q m c a_y u g) \, \dd \mu_{U / \G_U^g} (u) - \mu_{\Omega}(\phi) \Bigg | \, \dd \mu_{M}(m) \lesssim \| \phi \|_{L^\infty(\Omega)} \, \delta^{c_1}.
\end{equation}
Denote by $\mathfrak{S}(\delta)^c$ the complementary subset of $\mathfrak{S}(\delta)$ in $\mathfrak{S}$. By \eqref{eq:Siegel-Adjoint-Bound}, for every $m \in \mathfrak{S}(\delta)^c$, we have $\|\Ad(m) \| \lesssim \delta^{-1}$. Moreover, by \cite[Section 3.7]{EMV09}, provided $\ell$ is large enough, there is $c_2 > 0$ such that, for every $g \in G$ and $\phi \in C_c^\infty(\Omega)$, we have $\cS_{\ell} (g \cdot \phi) \lesssim \|\Ad(g) \|^{c_2} \, \cS_{\ell}(\phi)$ and $\|\phi \|_{L^\infty(\Omega)} \lesssim \cS_{\ell}(\phi)$. Applying \eqref{eq:Equi4} with the compact subset $Q' = Q \, \mathfrak{S}(\delta)^c \, C$ of $G$ (but keeping the dependency on $m \in \mathfrak{S}(\delta)^c$ explicit), we have
\begin{align} \label{align:Proof-Equi-Margulis}
\int_{\mathfrak{S}(\delta)^c} \sum_{c \in C} &\Bigg | \int_{U / \G_U^g} \phi(q m c a_y u g) \, \dd \mu_{U / \G_U^g} (u) - \mu_{\Omega}(\phi) \Bigg | \, \dd \mu_{m}(m) \notag \\ \notag
&\lesssim_{g, G, \G, Q} \int_{\mathfrak{S}(\delta)^c} \sum_{c \in C}  y^{c'} \, \|\Ad(m) \|^{c_2} \mathcal{S}_{\ell}(\phi) \, \dd \mu_{M}(m) \\ \notag
&\lesssim_{g, G, \G, Q} \int_{\mathfrak{S}(\delta)^c} y^{c'} \, \delta^{- c_2} \mathcal{S}_{\ell}(\phi) \, \dd \mu_{M}(m) \\ 
&\lesssim_{g, G, \G, Q} y^{c'} \, \delta^{- c_2} \mathcal{S}_{\ell}(\phi)
\end{align}
Therefore, combining \eqref{eq:Proof-Equi-Margulis} and \eqref{align:Proof-Equi-Margulis}, we have
\[
\Bigg | \int_{L / \G_L^g} \phi(k a_y l) \, \dd \mu_{L / \G_L^g} (l) - \mu_{\Omega}(\phi)  \Bigg | \lesssim y^{c'} \, \delta^{- c_2} \mathcal{S}_{\ell}(\phi) + \delta^{c_1} \mathcal{S}_{\ell}(\phi).
\]
Setting $\delta = y^{\frac{c'}{c_1 + c_2}}$ proves Lemma \ref{lem:Equi2} with $c = c' \left ( 1 - \frac{c_2}{c_1 + c_2} \right )$.
\end{proof}

\subsection{A lattice point counting lemma}
The following lemma is one of the key ingredients to Theorem \ref{thm:ThmMain2}.

\begin{lemma} \label{lem:Asymptotic}
There exist constants $\varkappa_2 > 0$, $\varepsilon_2 > 0$, and $T_0 \geq 1$ such that for all $T \geq T_0$ and all $g_1 \in G$ with $\mathrm{ht}(g_1\G) < T^{\varepsilon_2}$,
\[
\# \left (g_1 \cL_{\chi} \cap \cB_T \right ) = \varkappa_2 \, m_{\cV_{\chi}}(\cB_T) \left (1 + O( T^{-\varepsilon_2}) \right ).
\]
\end{lemma}
\begin{proof} 
We analyze the contribution of each $\G$-orbit in $\cL_{\chi}$ separately. In fact, by a theorem of Borel and Harish-Chandra \cite[Proposition~15.6]{Borel69}, the set of double cosets $\G \backslash \bG(\Q) / \bP(\Q)$ is finite. Moreover, according to \cite[Lemma 2.6]{BT65}, one has $(\bG / \bP)(\Q) = \bG(\Q) / \bP(\Q)$. As a consequence, $\bX(\Q)$ is a finite union of $\G$-orbits and, since there is a one-to-one correspondence between $\bX(\Q)$ and lines passing through elements of $\cL_{\chi}$, there exist finitely many $\bm{v}_1, \dots, \bm{v}_\kappa \in \cL_{\chi}$ such that
\begin{equation} \label{eq:GammaOrbits}
\cL_{\chi} = \bigsqcup_{i=1}^\kappa \G \bm{v}_i,
\end{equation}
and we can pick $\tau_i \in \bG(\Q)$ and $\lambda_i > 0$ such that $\bm{v}_i = \lambda_i \tau_i \bm{e}_\chi$. Note that $L_i = \tau_i L \tau_i^{-1}$ is the stabilizer of $\bm{v}_i$ in $G$ and put $\G_{L_i} = \G \cap L_i$. Then, for every $g \in G$, we have
\begin{equation}
\# (g \cL_{\chi} \cap \cB_T) = \sum_{i=1}^\kappa \# (g \G \bm{v}_i \cap \cB_T) = \sum_{i=1}^\kappa \sum_{\gamma \in \G/{\G_{L_i}}} \mathbbm{1}_{\cB_T} ( g \gamma \bm{v}_i).
\end{equation}
We now fix $1 \leq i \leq \kappa$ and define the function $F_T^{(i)} : \Omega \rightarrow \R$ by
\[
F_T^{(i)}(g\G) = \sum_{\gamma \in \G / \G_{L_i}} \mathbbm{1}_{\cB_T} (g \gamma  \bm{v}_i).
\]
Using Weil's integration formula \cite[Section 9]{Weil65} and the change of variable $g \tau_i \mapsto g$, for every measurable bounded function $\phi : \Omega \rightarrow \R$, we have
\begin{align*}
\langle F_T^{(i)}, \phi \rangle_{L^2(\Omega)} &= \int_{G / \G} \sum_{\gamma \in \G / \G_{L_i} } \mathbbm{1}_{\cB_T} (  g \gamma  \bm{v}_i) \phi(g \G) \, \dd  \mu_{\Omega} (g\G) \\ 
&= \int_{G / \G_{L_i}} \mathbbm{1}_{\cB_T} ( g \bm{v}_i) \phi(g \G) \, \dd \mu_{G/\G_{L_i}}( g \G_{L_i})\\
&= \int_{G / \tau_i^{-1} \G_{L_i} \tau_i} \mathbbm{1}_{\cB_T} ( \lambda_i g \bm{e}_\chi) \phi(g \tau_i^{-1} \G) \, \dd  \mu_{G/\tau_i^{-1} \G_{L_i} \tau_i} (g (\tau_i^{-1} \G_{L_i} \tau_i)).
\end{align*}
Let us write $\G_i = \tau_i^{-1} \G \tau_i$ and hence $\tau_i^{-1} \G_{L_i} \tau_i = \G_i \cap L$. Using Weil's integration formula again, together with the measure decomposition of $\mu_G$ and the fact that $L$ fixes $\bm{e}_\chi$, we have
\[
\langle F_T^{(i)}, \phi \rangle_{L^2(\Omega)} = \int_{\cV_{\chi}} \int_{L / \G_{i} \cap L} \mathbbm{1}_{\cB_T} ( \lambda_i k a_y \bm{e}_\chi) \phi(k a_y \ell \tau_i^{-1}) \, \dd \mu_{L / \G_{i} \cap L} (\ell \G_{i} \cap L) \dd m_{\cV_{\chi}}(k a_y \bm{e}_{\chi}).
\]
Since $\lambda_i \bm{e}_\chi = y_i^{-\frac{1}{\beta_\chi}} \bm{e}_\chi = a_{y_i} \bm{e}_\chi$ with $y_i = \lambda_i^{-\beta_\chi}$, by the change of variable $yy_i \mapsto y$, we have that $\langle F_T^{(i)}, \phi \rangle_{L^2(\Omega)}$ is given by
\[
\int_{\cV_{\chi}} \mathbbm{1}_{\cB_T} (  k a_y \bm{e}_\chi) \left ( \lambda_i^{-\beta_\chi d} \int_{L / \G_{i} \cap L} \phi(k a_{y_i}^{-1} a_y \ell \tau_i^{-1}) \, \dd \mu_{L / \G_{i} \cap L} (\ell) \right ) \dd m_{\cV_{\chi}}(k a_y \bm{e}_{\chi}).
\]
In particular, setting $\phi = 1$, we have  
\begin{equation} \label{eq:Mean-Value-Formula}
\int_{\Omega} F_T^{(i)} \, \dd \mu_{\Omega} = \nu_i \, m_{\cV_{\chi}}(\cB_T), \quad \text{with } \, \nu_i = \lambda_i^{-\beta_\chi d} \, \mu_{L/(\G_i \cap L)} \left (L/(\G_i \cap L) \right ).
\end{equation}
Hence we only need to show that there exists $\varepsilon_2 > 0$, independent of $i$, such that for all sufficiently large $T$ and for all $g_1 \in G$ with $\mathrm{ht}(g_1\G) < T^{\varepsilon_2}$, we have
\[
F_T^{(i)}(g_1\G) = \# \left (  g_1 \G  \bm{v}_i \cap \cB_T \right ) = \nu_i \, m_{\cV_{\chi}}(\cB_T) \left ( 1 + O( T^{-\varepsilon_2}) \right ).
\]
Let $\omega_i = \lambda_i^{\beta_\chi d} \nu_i$ and we will obtain Lemma \ref{lem:Asymptotic} with 
\begin{equation} \label{eq:Constant-In-Lemma-Asymptotic}
\varkappa_2 = \sum_{i = 1}^\kappa \nu_i = \sum_{i = 1}^\kappa \lambda_i^{-\beta_\chi d} \mu_{L/(\G_i \cap L)} \left (L/(\G_i \cap L) \right ). 
\end{equation}
Note that
\[
\sup \, \{ y \in \R_+^\times: \exists \, k \in K, \, k a_y \bm{e}_\chi \in \cB_T \} \lesssim T^{-(\beta_\chi -\tau)}.
\]
By Lemma \ref{lem:Equi2} (applied with $g = \tau_i^{-1}$ and the compact set $Q = K$), there exists $c > 0$ such that, for all $\phi \in C_c^\infty(\Omega)$, we have
\begin{align*}
\Bigg | &\langle F_T^{(i)}, \phi \rangle_{L^2(\Omega)} - \nu_i \, m_{\cV_{\chi}}(\cB_T) \int_{\Omega} \phi \, \dd \mu_{\Omega} \Bigg | \\
&\leq \int_{\cV_{\chi}} \mathbbm{1}_{\cB_T} (  k a_y \bm{e}_\chi) \lambda_i^{- \beta_\chi d} \, \Bigg | \int_{L / \G_{i} \cap L} \phi(k a_{y_i}^{-1} a_y \ell \tau_i^{-1}) \, \dd \mu_{L / \G_{i} \cap L} (\ell) \\
& \quad \quad \quad \quad \quad \quad \quad \quad \quad \quad \quad \quad \quad \quad \quad \quad \quad \quad - \omega_i \int_{\Omega} \phi \, \dd \mu_{\Omega} \Bigg | \, \dd m_{\cV_{\chi}} (k a_y \bm{e}_{\chi}) \\
&\lesssim \left (\int_{\cV_{\chi}} \mathbbm{1}_{\cB_T} (  k a_y \bm{e}_\chi) y^{c} \, \dd m_{\cV_{\chi}}(k a_y \bm{e}_{\chi})\right ) \, \cS_{\ell}(\phi) \\ 
&\lesssim T^{-c (\beta_\chi -\tau)} \, m_{\cV_{\chi}} (\cB_T) \, \cS_{\ell}(\phi).
\end{align*}
Let $c_1 > 0$ be as in Lemma \ref{lem:WellRounded}, and let $\cO = (\cO_\delta)_{0<\delta<1}$ with
\begin{equation*} 
\cO_\delta = B_G(c_1 \delta).
\end{equation*}    
We recall from Section \ref{sec:Sobolev} that, for every $x \in \Omega$, $c_3 \mathrm{ht}(x)^{-\kappa_0}$ is a lower bound for the injectivity radius $\mathrm{inj}(x)$ at $x$. Let $\varepsilon_2 > 0$ and $g_1 \in G$ with $\mathrm{ht}(g_1 \G) < T^{\varepsilon_2}$. By Lemma \ref{lem:Bump-Function-Cusp}, applied with $x = g_1 \G$ and $\delta = c_3 T^{-\varepsilon_2 \kappa_0} \leq c_3 \mathrm{ht}(g_1 \G)^{-\kappa_0} \leq \mathrm{inj}(g_1 \G)$, there exists a non-negative function $\phi \in C_c^\infty(\Omega)$ such that
\[
\mu_{\Omega}(\phi) = 1, \quad \supp(\phi) \subseteq \cO_\delta \, g_1 \G, \quad \text{and} \quad \cS_{\ell} (\phi) \lesssim \delta^{-D} \mathrm{ht}(x)^{\ell} \lesssim T^{-\varepsilon_2 (\kappa_0 (\ell + \dim G / 2) + \ell)}
\]
Putting everything together and letting $D' = \kappa_0 (\ell + \dim G / 2) + \ell$, we have
\[
\left | \langle F_T^{(i)}, \phi \rangle_{L^2(\Omega)} - \nu_i \,  m_{\cV_{\chi}} (\cB_T) \right | \lesssim T^{-c (\beta_\chi -\tau) + \varepsilon_2 D'} m_{\cV_{\chi}} (\cB_T).
\]
Using Proposition \ref{prop:WellRounded}, the family $(\cB_T)_{T \geq 1}$ is well-rounded, and we showed that there exists $C_1 > 0$ such that for all sufficiently large $T \geq 1$ and $\delta \in (0,1)$, there are Borel sets $\underline{\cB}_{T,\delta}, \overline{\cB}_{T,\delta}$ satisfying
\begin{equation*} 
\underline{\cB}_{T,\delta}\subseteq \bigcap_{g\in \cO_{\delta}}g\cB_T \subseteq \bigcup_{g\in \cO_{\delta}}g \cB_T \subseteq \overline{\cB}_{T,\delta}
\end{equation*}
and $m_{\cV_{\chi}}(\overline{\cB}_{T,\delta}\smallsetminus \, \underline{\cB}_{T,\delta}) \leq C_1 \delta m_{\cV_{\chi}}(\cB_T)$. For $\delta = c_3 T^{- \varepsilon_2 \kappa_0}$ define 
\[
\underline{F}_{T,\delta}(g) = \sum_{\gamma \in \G / \G_{L_i}} \mathbbm{1}_{\underline{\cB}_{T,\delta}} (  g \gamma  \tau_i \bm{e}_\chi)
\quad \text{and} \quad \overline{F}_{T,\delta}(g) = \sum_{\gamma \in \G / \G_{L_i}} \mathbbm{1}_{\overline{\cB}_{T,\delta}} (  g \gamma  \tau_i \bm{e}_\chi).
\]
Then for $g\G \in \supp \, \phi$, there exists $h \in \cO_\delta$ such that  $g \G = h g_1 \G$. Thus
\[
\underline{F}_{T,\delta}(g) = \underline{F}_{T,\delta}(h g_1) = 
\# \left (  g_1 \G  \tau_i \bm{e}_\chi \cap h^{-1} \underline{\cB}_{T,\delta} \right ) \leq F_T^{(i)}( g_1).
\]
Multiplying by $\phi(g)$ and integrating over the support of $\phi$ gives
\begin{equation*}
\langle \underline{F}_{T,\delta}, \phi \rangle_{L^2(\Omega)} \leq F_T^{(i)}(g_1).
\end{equation*}
Similarly as above, we can estimate
\[
\left | \langle \underline{F}_{T,\delta}, \phi \rangle_{L^2(\Omega)} - \nu_i m_{\cV_{\chi}}(\underline{\cB}_{T,\delta}) \right | \lesssim T^{-c (\beta_\chi -\tau) + \varepsilon_2 D'}  m_{\cV_{\chi}}(\underline{\cB}_{T,\delta}).
\]
Moreover, from $m_{\cV_{\chi}}(\overline{\cB}_{T,\delta}\smallsetminus \,\underline{\cB}_{T,\delta}) \leq C_1 \delta m_{\cV_{\chi}}(\cB_T)$ we get 
\[
m_{\cV_{\chi}}(\underline{\cB}_{T,\delta})\geq (1 - C_1 \delta) m_{\cV_{\chi}}(\cB_{T}).
\]
Therefore, one has
\begin{align*}
F_T^{(i)}(g_1) - \nu_i m_{\cV_{\chi}}(\cB_{T}) &\geq \langle \underline{F}_{T,\delta}, \phi \rangle_{L^2(\Omega)} - \nu_i m_{\cV_{\chi}}(\underline{\cB}_{T,\delta}) + \nu_i m_{\cV_{\chi}}(\underline{\cB}_{T,\delta}) - \nu_i m_{\cV_{\chi}}(\cB_{T}) \\
&\gtrsim - T^{-c (\beta_\chi -\tau) +\varepsilon_2 D'}  m_{\cV_{\chi}}(\underline{\cB}_{T,\delta}) - C_1 \delta m_{\cV_{\chi}}(\cB_{T}) \\
&\gtrsim -(T^{-c (\beta_\chi -\tau) +\varepsilon_2 D'} + T^{-\varepsilon_2 \kappa_0}) m_{\cV_{\chi}}(\cB_{T}). 
\end{align*}
Similarly, using the subset $\overline{\cB}_{T,\delta}$, one can show that 
\[
F_T^{(i)}(g_1) - \nu_i m_{\cV_{\chi}}(\cB_{T}) \lesssim (T^{-c (\beta_\chi -\tau) +\varepsilon_2 D'} + T^{-\varepsilon_2 \kappa_0}) m_{\cV_{\chi}}(\cB_{T}). 
\]
Together, this gives 
\[
\left | F_T^{(i)}(g_1) - \nu_i m_{\cV_{\chi}}(\cB_{T}) \right | \lesssim (T^{-c (\beta_\chi -\tau) +\varepsilon_2 D'} + T^{-\varepsilon_2 \kappa_0} ) m_{\cV_{\chi}}(\cB_{T}).
\]
This proves Lemma \ref{lem:Asymptotic} with $\varepsilon_2 = \frac{c (\beta_\chi -\tau)}{D'+\kappa_0}$ and $\varkappa_2$ given by Equation \eqref{eq:Constant-In-Lemma-Asymptotic}.
\end{proof}

\subsection{Proof of Proposition \ref{prop:Counting}}
\begin{proof} [Proof of Proposition \ref{prop:Counting}]
By the proof of \cite[Theorem~2.4.5]{deSaxce20}, using the ergodicity of the action of $A$ on the probability space $(\Omega, \mu_{\Omega})$, for $\mu_G$-almost every $g \in G$, we have
\begin{equation} \label{eq:Limit-Ergodic}
\lim_{y \rightarrow \infty} \frac{\ln \left (\lambda_1(  \Ad(a_y g) \kg_\G ) \right )}{\ln(y)}  = 0.
\end{equation}
For every $p \in P$, there exists $p_\infty \in P$ such that, as $y \rightarrow + \infty$, we have $a_y p a_y^{-1} \rightarrow p_\infty$. Hence the above limit only depends on the class $Pg$ of $g$ in $P \backslash G$. Therefore, for $\sigma_X$-almost every $x = k_x P \in X$, we have
\begin{equation} \label{eq:Limit-Ergodic-2}
\lim_{y \rightarrow \infty} \frac{\ln \left (\lambda_1(  \Ad(a_y k_x^{-1}) \kg_\G ) \right )}{\ln(y)}  = 0.
\end{equation}
Plugging in the definition $y_T = T^\tau$ of $y_T$, there exists a full measure subset $X_0 \subseteq X$ such that for every $x \in X_0$, we have, as $T \rightarrow + \infty$,
\begin{equation} \label{eq:Ergodic_Relation}
\lambda_1(  \Ad(a_{y_T} k_x^{-1}) \kg_\G ) = T^{o(1)}.
\end{equation}
Fix $x \in X_0$. Since $\cL_{\chi} = \bigsqcup_{1 \leq i \leq \kappa} \G \, \bm{v}_i$ is a finite union of $\G$-orbits (see \eqref{eq:GammaOrbits}), it is sufficient to prove that there exists $\varepsilon \in (0,1)$, independent of $i$, such that for all sufficiently large $T \geq 1$,
\begin{equation} \label{eq:Final-Equation}
\# \left (  k_x^{-1} \G \bm{v}_i \cap \cE_\tau(T) \right ) =  \nu_i \, m_{\cV_{\chi}}(\cE_\tau(T)) \left (1 + O_x ( T^{-\varepsilon} ) \right),
\end{equation}
where $\nu_i$ is as in Equation \eqref{eq:Constant-In-Lemma-Asymptotic}.
Recall that for every $T \geq 1$, we defined
\[
\cF_T = \left \{ \bm{v} \in \cV_{\chi}: d(x_0, [\bm{v}]) < \psi(\|\bm{v}\|), \, \max \{1 ,T/2 \} \leq \|\bm{v}\| < T \right \}.
\]
This gives the disjoint union
\[
\cE_\tau(T) = \bigsqcup_{j \geq 0} \cF_{T_j}, \quad \text{ with } T_j = T /2^{j} \text{ for every  $j \in \N$}, 
\]
and therefore 
\begin{equation} \label{eq:1_prop:Counting}
\# \left (  k_x^{-1} \G \bm{v}_i \cap \cE_\tau(T) \right ) = \sum_{j \geq 0 } \# \left (  k_x^{-1} \G \bm{v}_i \cap \cF_{T_j} \right ).
\end{equation}
For a given $T \geq 1$, we now analyze and estimate each term $\# \left (  k_x^{-1} \G \bm{v}_i \cap \cF_{T} \right )$ individually. In view of \eqref{eq:Ergodic_Relation} we have $\mathrm{ht}(a_{y_T} k_x^{-1} \G) = T^{o(1)}$. Since $\cB_T = a_{y_T} \cF_T$ and the measure $m_{\cV_{\chi}}$ is invariant under $A$, we have $m_{\cV_{\chi}}(\cB_T) = m_{\cV_{\chi}}(\cF_T)$. Now by Lemma \ref{lem:Asymptotic}, there exists $\varepsilon_2 > 0$ such that for all sufficiently large $T > T_0$,
\begin{align*}
\# \left (  k_x^{-1} \G \bm{v}_i \cap \cF_{T} \right ) &= \# \left (  a_{y_T} k_x^{-1} \G \bm{v}_i \cap \cB_{T} \right ) \\
&= \nu_i \, m_{\cV_{\chi}}(\cF_T)\left ( 1 + O_x( T^{-\varepsilon_2} ) \right ) \\
&= \nu_i \, m_{\cV_{\chi}}(\cB_T)\left ( 1 + O_x( T^{-\varepsilon_2} ) \right ).
\end{align*}
Note that
\[
\sum_{T_j \leq T_0} \# \left (  k_x^{-1} \G \tau_i \bm{e}_\chi \cap \cF_{T_j} \right ) - \sum_{T_j \leq T_0} \nu_i \, m_{\cV_{\chi}}(\cF_{T_j}) = O_x\left (\sum_{j \geq 0} T_j^{-\varepsilon_2} m_{\cV_{\chi}} (\cF_{T_j}) \right ).
\]
Returning to \eqref{eq:1_prop:Counting}, we deduce that
\begin{align*}
\# &\left (  k_x^{-1} \G \bm{v}_i \cap \cE_\tau(T) \right ) = \sum_{j \geq 0 } \# \left (  k_x^{-1} \G \bm{v}_i \cap \cF_{T_j} \right ) \\
&= \sum_{T_j \leq T_0} \# \left (  k_x^{-1} \G \bm{v}_i \cap \cF_{T_j} \right ) + \sum_{T_j > T_0} \nu_i \, m_{\cV_{\chi}}(\cF_{T_j}) + O_x \left (\sum_{T_j > T_0} T_j^{-\varepsilon_2} m_{\cV_{\chi}} (\cF_{T_j}) \right ) \\
&= \nu_i \, m_{\cV_{\chi}}(\cE_\tau(T)) +O_x\left (\sum_{j \geq 0} T_j^{-\varepsilon_2} m_{\cV_{\chi}} (\cF_{T_j}) \right ),
\end{align*}
Thus to prove the formula \eqref{eq:Final-Equation}, it suffices to find $\varepsilon_1 \in (0,1)$ such that
\[
\sum_{j \geq 0} T_j^{-\varepsilon_2} m_{\cV_{\chi}} (\cF_{T_j}) \lesssim T^{-\varepsilon_1} m_{\cV_{\chi}} (\cE_\tau(T)).
\]
Using the definition $T_j = T / 2^j$ and the volume estimate 
\[
m_{\cV_{\chi}} (\cF_{T_j}) \lesssim \int_{T_{j+1}}^{T_j} y^{(\beta_\chi -\tau) d} \frac{d y}{y} \lesssim T_j^{(\beta_\chi - \tau) d} - T_{j+1}^{(\beta_\chi - \tau) d} = (1 - 2^{-(\beta_\chi - \tau) d})T_j^{(\beta_\chi - \tau) d},
\] 
and setting $\varepsilon_1 = \min\{\varepsilon_2, (\beta_\chi - \tau)d /2\}$, $\sum_{j \geq 0} T_j^{-\varepsilon_2} m_{\cV_{\chi}} (\cF_{T_j})$ is bounded by
\begin{align*}
\sum_{j \geq 0} T_j^{(\beta_\chi - \tau) d - \varepsilon_2} \lesssim T^{(\beta_\chi - \tau) d -\varepsilon_1} \sum_{j \geq 0} \left ( \frac{1}{2^{(\beta_\chi - \tau)d-\varepsilon_1}} \right )^j
\lesssim T^{-\varepsilon_1} m_{\cV_{\chi}}(\cE_\tau(T)).
\end{align*}
This establishes the proof of Proposition \ref{prop:Counting} with $\varepsilon_1 = \min\{\varepsilon_2, (\beta_\chi - \tau)d /2\}$ and $\varkappa_1 = \varkappa_2$ as in Lemma \ref{lem:Asymptotic}.
\end{proof}

\section{Counting at the Diophantine exponent} \label{sec:Critical}
In this section, we prove Theorem \ref{thm:ThmMain}, where we count rational approximations at the Diophantine exponent $\beta_{\chi}$ of $X$ relative to $\chi$ (see Section \eqref{sec:Diophantine}). 

\medskip
The method is inspired by the ergodic-theoretic approach in \cite{AG22}, where Alam and Ghosh counted rational approximations on spheres. Let us briefly describe this method.
We recall that, by Lemma \ref{lem:Reduction}, for every $x \in X$ and $T \geq 1$, we have $\cN_{c,\beta_\chi}(x,T) = \# \left (k_x^{-1} \cL_{\chi} \cap \cE_{c, \beta_{\chi}}(T) \right )$, where 
\[
\cE_{c, \beta_{\chi}}(T) = \{\bm{v} \in \cV_{\chi} : d(x_0, [\bm{v}]) < c \, \|\bm{v}\|^{-\beta_\chi}, 1 \leq \|\bm{v}\| < T\}.
\]
First, one approximates $\cE_{c, \beta_{\chi}}(T)$ by regions that admit a decomposition with respect to the action of the diagonal subgroup $A$. To conclude, one uses Birkhoff’s ergodic theorem and an approximation argument. We assume for simplicity that $c = 1$ and put $\cE_{\beta_{\chi}}(T) = \cE_{1, \beta_{\chi}}(T)$; for general $c > 0$ the argument is identical.

\medskip
Let $\pi^+ : V_\chi \rightarrow V_\chi$ be the orthogonal projection onto $\R\bm{e}_\chi$ and we abbreviate $\pi^+(\bm{v})$ simply by $\bm{v}^+$. Recall from \textsection \ref{sec:CC} that $\phi : G \rightarrow X$ is the projection map sending $g \mapsto g x_0$. Denote by $D_1 \phi : \kg \rightarrow T_{x_0} X$ its derivative at the identity $1 \in G$. It satisfies $\ker D_1 \phi = \mathfrak{p}$, where $\mathfrak{p}$ is the Lie algebra of $P$, and hence defines an isomorphism $D_1 \phi : \ku^- \rightarrow T_{x_0} X$. We use this isomorphism to identify $\ku^-$ with $T_{x_0} X$. By \cite[Lemme~2.4.2]{deSaxce20}, using the definition 
\[
\forall \, y > 0, \quad a_{y} = \exp(\ln (y) Y_\alpha)
\]
and the fact that $\ku^-$ is abelian, the adjoint action of $a_{y}$ on $\ku^- = T_{x_0} X$ is given by
\[
\forall \, u \in \ku^-, \quad \Ad(a_{y}) u^- = y u^-.
\]
For every non-zero $\bm{v} \in \cV_{\chi}$, such that $[\bm{v}]$ is close to $x_0$, we denote by $u_{\bm{v}}^- \in \ku^-$ the unique element such that $[\bm{v}] = \exp(u_{\bm{v}}^-) x_0$. Observe that 
\[
[a_y \bm{v}] = a_y [\bm{v}] = a_y \exp(u_{\bm{v}}^-) a_y^{-1} a_y x_0 = \exp(\Ad(a_y) u_{\bm{v}}^-) x_0 = \exp(y u_{\bm{v}}^-) x_0
\]
But one also has $[a_y \bm{v}] = \exp(u_{a_y \bm{v}}^- ) x_0$. By uniqueness, this gives the relationship
\begin{equation} \label{eq:LieAlgebraDiagonalErgodic}
u_{a_y \bm{v}}^- = y \, u_{\bm{v}}^-.
\end{equation}
By \eqref{eq:Distance_Estimate}, there is a constant $C_0 > 0$ such that 
\begin{equation} \label{eq:Distances}
d(x_0, [\bm{v}]) \leq \|u_{\bm{v}}^-\|_{\ku^-} + C_0 \|u_{\bm{v}}^-\|_{\ku^-}^2.
\end{equation}
We now approximate the region $\cE_{\beta_{\chi}}(T)$ from inside and from outside by regions
\begin{equation} \label{eq:New_Region}
\cE_{T, c}^+ = \{\bm{v} \in \cV_{\chi} : \|u_{\bm{v}}^-\|_{\ku^-} < c \, \|\bm{v}^+\|^{-\beta_\chi}, 1 \leq \|\bm{v}^+\| < c \, T \},
\end{equation}
where $c>0$ is a parameter that will approach $1$. By enlarging $C_0$ if necessary, we can assume that $\|\bm{v}^+\| \geq C_0^{-1} \|\bm{v}\|$ if $[\bm{v}]$ is sufficiently close to $x_0$. For every natural number $\ell \geq 1$, let 
\begin{equation} \label{eq:Defintion-Q_l}
Q_\ell = \{\bm{v} \in \cV_{\chi} : \|\bm{v}\| \leq C_0 \ell \} \,\, \text{ and } \,\, c_\ell = \left ( 1 + C_0 \, \ell^{-\beta_\chi/2}  \right )^{-2(1+\beta_\chi)} \in (0,1).
\end{equation}
In particular, one has $c_\ell \nearrow 1$ as $\ell \rightarrow + \infty$. Fix a right-invariant Riemannian metric $d_G(\cdot, \cdot)$ on $G$. Write $B_P(r)$ for the (symmetric) open ball in $P$ with radius $r > 0$ and center $1 \in P$ with respect to the distance induced by that on $G$. 

\medskip
The following lemma tells us that the region $\cE_{\beta_{\chi}}(T)$ can be approximated well by regions of the form $\cE_{T,c}^{+}$, that this approximation is even stable under perturbations of elements $p$ close to the identity in $P$, and that the regions $\cE_{T,c}^{+}$ admit a decomposition with respect to the action of the subgroup $A$. For every $c > 0$, let 
\begin{equation} \label{eq:F_c}
\cF_c = \{\bm{v} \in \cV_{\chi} : \|u_{\bm{v}}^-\|_{\ku^-} < c \, \|\bm{v}^+\|^{-\beta_{\chi}}, \, 1 \leq \|\bm{v}^+\| < 2 \}.
\end{equation}

\begin{lemma} \label{lem:SandwichErgodic}
For all large enough $\ell \geq 1$, $p \in \cO_\ell = B_P(\ell^{-\beta_\chi/2})$, and $T \geq 1$,
\begin{equation} \label{eq:lem:SandwichErgodic}
\cE_{T, c_\ell}^+ \setminus Q_{2\ell} \, \subseteq \, p \, \left ( \cE_{\beta_{\chi}}(T) \setminus Q_{\ell} \right ) \, \subseteq \, \cE_{T, c_\ell^{-1}}^+.
\end{equation}
Moreover, for all $c > 0$ and $T \geq 1$ such that $c T = 2^N$ for some $N \in \N$, we have
\begin{equation} \label{eq:TessellationErgodic}
\cE_{T,c}^+ = \bigsqcup_{i=0}^{N-1} a_{y_j}^{-1}\, \cF_c, \quad \quad \text{with $y_j = 2^{\beta_{\chi} j}$ for every $j \in \N$. }
\end{equation}
\end{lemma}

\begin{proof}
Let us show the left inclusion in Equation \eqref{eq:lem:SandwichErgodic}. Put $\delta_\ell = \ell^{-\beta_\chi/2}$. We need to show that for all large enough $\ell \geq 1$, $p \in \cO_\ell$, $T \geq 1$, and $\bm{v} \in \cE_{ T, c_\ell}^+ \setminus Q_{2 \ell}$,
\[
d(x_0, p[\bm{v}]) < \|p\bm{v}\|^{-\beta_\chi}, \quad \text{and} \quad C_0 \ell < \| p\bm{v}\| < T.
\]
Using the triangle inequality, we get $d(x_0, p[\bm{v}]) \leq d(x_0, [\bm{v}]) + d([\bm{v}], p[\bm{v}])$. Next, writing $[\bm{v}] = \exp(u_{\bm{v}}^-) x_0$ with $u_{\bm{v}}^- \in \ku^-$, using that $p$ stabilizes the line $x_0$, and expressing $\exp(\Ad(p)u_{\bm{v}}^-) = \exp(u')p'$ with $u' \in \ku^-, \, p' \in P$, we get
\begin{align*}
d([\bm{v}], p[\bm{v}]) &\lesssim  d(\exp(u_{\bm{v}}^-)x_0,  p \exp(u_{\bm{v}}^-) x_0) \\
&=  d(\exp(u_{\bm{v}}^-) x_0,  \exp(\Ad(p)u_{\bm{v}}^-) x_0) \\
&=  d(\exp( u_{\bm{v}}^-) x_0, \exp(u') x_0) \asymp \|u' - u_{\bm{v}}^-\|_{\ku^-}
\end{align*}
Note that the map from a small neighborhood of $0$ in the Lie algebra $\kg$ of $G$ to $\ku^-$ defined by $X \mapsto X'$ with $X' \in \ku^-$ and $p_X \in P^\circ$ such that $\exp(X) = \exp(X') p_X$ is Lipschitz. Hence
\begin{multline*}
\|u' - u_{\bm{v}}^-\|_{\ku^-} \lesssim \| \Ad(p)u_{\bm{v}}^- - u_{\bm{v}}^-\|_{\kg} \lesssim \|\Ad(p) - \mathrm{Id} \| \|u_{\bm{v}}^-\|_{\ku^-} \lesssim \delta_\ell \, d(x_0, [\bm{v}]).
\end{multline*}
Therefore, by enlarging $C_0$ if necessary, one has $d(x_0, p[\bm{v}]) \leq (1 + C_0 \delta_\ell) d(x_0, [\bm{v}])$.  Together with the estimate \eqref{eq:Distances} and the fact that $\bm{v} \in \cE_{ T, c_\ell}^+$, we thus get
\begin{align*} 
d(x_0, p[\bm{v}]) &\leq (1 + C_0 \delta_\ell) \, \|u_{\bm{v}}^-\|_{\ku^-} \, \left ( 1 + C_0 \|u_{\bm{v}}^-\|_{\ku^-} \right ) \\
&< \|\bm{v}^+ \|^{-\beta_\chi} c_\ell \, (1 + C_0 \delta_\ell)  \left ( 1 + C_0 \|\bm{v}^+ \|^{-\beta_\chi} \right ).
\end{align*}
Let $\bm{v}^{\perp} = \bm{v} - \bm{v}^+$ be the projection of $\bm{v}$ onto the orthogonal complement of $\R \bm{e}_{\chi}$. Using the fact that $\bm{v} \in \cE_{ T, c_\ell}^+ \setminus Q_{2 \ell}$, we have 
\[
\frac{\|\bm{v}^{\perp}\|}{\|\bm{v}^+ \|} \asymp d(x_0, v) \asymp \|u_{\bm{v}}^-\|_{\ku^-} \lesssim \|\bm{v}^+ \|^{-\beta_\chi}.
\]
Hence, by enlarging $C_0$ if necessary and using that $\|\bm{v}\|^2 = \|\bm{v}^+\|^2 + \|\bm{v}^{\perp}\|^2$, we have
\[
\|\bm{v}^+ \|^{-\beta_\chi} = \|\bm{v} \|^{-\beta_\chi} \left (1 + \frac{\|\bm{v}^{\perp} \|^2}{\|\bm{v}^+ \|^2} \right )^{\frac{\beta_\chi}{2}} \leq \|\bm{v} \|^{-\beta_\chi} \left (1 + C_0 \|\bm{v}^+ \|^{-2\beta_\chi} \right )^{\frac{\beta_\chi}{2}}.
\]
By enlarging $C_0$, we may assume that $\|\bm{v}^+\| \geq C_0^{-1} \|\bm{v}\|$. Since $\bm{v} \notin Q_{2\ell}$, we thus have $\|\bm{v}^+\| \geq C_0^{-1} \|\bm{v}\| \geq 2 \ell$. Moreover, we may assume that $C_0$ is such that $\|p \bm{v} \| \leq (1 + C_0 \delta_\ell) \|\bm{v} \|$ for all large $\ell \geq 1$ and $p \in B_P(\delta_\ell)$. Putting everything together and using the definition of $c_\ell$, we have
\begin{align*}
d(x_0, p[\bm{v}]) &< \|p \bm{v}\|^{-\beta_\chi} \left ( c_\ell (1 + C_0 \delta_\ell)^{1 + \beta_\chi} (1 + C_0 \|\bm{v}^+ \|^{-\beta_\chi}) \left (1 + C_0 \|\bm{v}^+ \|^{-2\beta_\chi} \right )^{\frac{\beta_\chi}{2}} \right ) \\
&\leq \|p \bm{v}\|^{-\beta_\chi} \left ( c_\ell (1 + C_0 \delta_\ell)^{1 + \beta_\chi} (1 + C_0 \ell^{-\beta_\chi}) \left (1 + C_0 \ell^{-2\beta_\chi} \right )^{\frac{\beta_\chi}{2}} \right ) \\
&\leq \|p \bm{v}\|^{-\beta_\chi},
\end{align*}
as desired. Moreover, we have 
\[
\|p \bm{v}\| \leq (1 + C_0 \, \delta_{\ell} ) \|\bm{v}^+\| \frac{\|\bm{v}\|}{\|\bm{v}^+\|} \leq  (1 + C_0 \, \delta_{\ell} ) c_\ell \frac{\|\bm{v}\|}{\|\bm{v}^+\|} T < T.
\]
Finally, using that $\bm{v} \notin Q_{2 \ell}$, we have, when $\ell$ is large,
\[
\|p \bm{v}\| \geq (1 + C_0 \delta_\ell)^{-1} \|\bm{v}\| > (1 + C_0 \delta_\ell)^{-1} C_0 (2 \ell) \geq C_0 \ell.
\]
This shows the left inclusion in Equation \eqref{eq:lem:SandwichErgodic}. The right inclusion is proved similarly. To see the last claim, let us recall that $a_y$ acts on $\bm{v}^+$ by $a_y \bm{v}^+ = y^{-\frac{1}{\beta_{\chi}}}\bm{v}^+$. Then the claim follows by using \eqref{eq:LieAlgebraDiagonalErgodic} and observing that
\[
a_{y_j}^{-1} \, \cF_c = \{\bm{v} \in \cV_{\chi} : \|u_{\bm{v}}^-\|_{\ku^-} < c \, \|\bm{v}^+\|^{-\beta_{\chi}}, 2^j \leq \|\bm{v}^+\| < 2^{j+1} \}. 
\]
The proof of Lemma \ref{lem:SandwichErgodic} is complete.
\end{proof}
We recall that, for every $j \in \Z$, we defined $y_j = 2^{\beta_{\chi} j}$. By Moore’s ergodicity theorem (see \cite[Section 3, Theorem 2.1]{BM00}), the action of the unbounded subgroup $\{a_{y_j} : j \in \Z \}$ of $G$ on the probability space $(\Omega, \mu_{\Omega})$ is ergodic. Therefore, by Birkhoff's ergodic theorem (see \cite[Section 1, Theorem 2.5]{BM00}), for every $f \in L^1(\Omega)$ and almost every $x \in \Omega$,
\begin{equation} \label{eq:Birkhoff}
\frac{1}{N} \sum_{j=0}^{N-1} f(a_{y_j} x) \longrightarrow \int_{\Omega} f \, \dd \mu_{\Omega} \quad \quad \text{ as $N \rightarrow + \infty$.}
\end{equation}
A point $x \in \Omega$ satisfying \eqref{eq:Birkhoff} is called \emph{Birkhoff generic} with respect to $f$. 

\begin{proof} [Proof of Theorem \ref{thm:ThmMain}]
For any $c > 0$ with $\cF_c$ as defined in Equation \eqref{eq:F_c}, define the function $F_c : \Omega \rightarrow \R$ by $F_c : g \G \mapsto \# (g \cL_{\chi} \cap \cF_{c})$. By the proof of Lemma \ref{lem:Asymptotic} and with $\varkappa_2 > 0$ as in \eqref{eq:Constant-In-Lemma-Asymptotic}, for every $c > 0$, we have
\begin{equation} \label{eq:First_Moment}
\int_{\Omega} F_c \, \dd \mu_{\Omega} =  \varkappa_2 \, m_{\cV_{\chi}}(\cF_{c}). 
\end{equation}
By the Iwasawa decomposition, the set $S = \cO_1 \cdot K$, where $\cO_1 = B_P(1)$ is as above, contains an open neighborhood of the identity in $G$. Moreover, a Fubini-type argument shows that for almost every $p \in \cO_1$ there is a measurable subset $K_p \subset K$ with $\mu_K(K_p) = 1$ such that for every $k \in K_p$ the point $p k \G$ is Birkhoff generic with respect to the function $F_c$. Consequently, for all $\ell \geq 1$, we can find $p_\ell \in \cO_\ell$ and a full-measure subset $K_\ell \subset K$ such that for every $k \in K_\ell$ the point $p_\ell k \G$ is Birkhoff generic with respect to the functions  $F_{c_\ell}$ and $F_{c_\ell^{-1}}$. Let $k \in K_\infty = \bigcap_{\ell \geq 1} K_\ell$. Since $p_\ell \in \cO_\ell$, Lemma \ref{lem:SandwichErgodic} gives for all large $\ell$ and $T \geq 1$, that
\[
\cE_{T, c_\ell}^+ \setminus Q_{2\ell} \, \subseteq \, p_\ell \cdot \left ( \cE_{\beta_{\chi}}(T) \setminus Q_{\ell} \right ) \, \subseteq \, \cE_{T, c_\ell^{-1}}^+.
\]
Intersecting with $p_\ell k \cL_{\chi}$ and using that the number of lattice points in the set $\cO_1 Q_\ell$ is bounded by an absolute constant times $\ell^{\beta_\chi d}$ (see Theorem \ref{thm:ThmMainEquid}), up to enlarging $C_0$ if necessary, we have
\begin{equation} \label{eq:Sandwich_lem:ReductionErgodic}
\# \left (p_\ell k \cL_{\chi}  \cap \cE_{T, c_\ell}^+ \right ) - C_0 \ell^{\beta_\chi d} \leq \# \left (k \cL_{\chi}  \cap \cE_{T} \right ) \leq \# \left (p_\ell k \cL_{\chi}  \cap \cE_{T, c_\ell^{-1}}^+ \right ) + C_0 \ell^{\beta_\chi d}.
\end{equation}
Using the decomposition of $\cE_{T,c}^+$ in \eqref{eq:TessellationErgodic} with $\ell$ large enough and $T = \frac{1}{c}2^N$ for every integer $N \geq 1$,  we have
\begin{equation} \label{eq:Method_Birkhoff}
\# (p_\ell k \cL_{\chi} \cap \cE_{T,c}^+) = \sum_{j=0}^{N-1} F_c (a_{y_j} p_\ell k \G).
\end{equation}
For $T \geq 1$ and $\ell$ large, let $N \geq 1$ be the unique integer with $\frac{1}{c_{\ell}} 2^{N} \leq T < \frac{1}{c_{\ell}} 2^{N+1}$. Plugging this back into \eqref{eq:Sandwich_lem:ReductionErgodic}, we get the lower and upper bounds
\begin{align*} \label{eq:Sandwich_lem:Reduction_Birkhoof}
\sum_{j=0}^{N-1} F_{c_\ell} (a_{y_j} p_\ell k \G) - C_0 \ell^{\beta_\chi d} &\leq \# \left ( k \cL_{\chi}  \cap \cE_{T} \right ) \notag \\ 
&\leq \sum_{j=0}^{N} F_{c_\ell^{-1}} (a_{y_j} p_\ell k \G) + C_0 \ell^{\beta_\chi d}.
\end{align*}
Dividing by $\ln (T)$, using that $p_\ell k \G$ is Birkhoff generic with respect to $F_{c_\ell}$ and $F_{c_\ell^{-1}}$, taking limits $T \rightarrow +\infty$ and evaluating them using \eqref{eq:First_Moment}, we have
\begin{equation} \label{eq:Sandwich_lem:Reduction2}
\frac{\varkappa_2}{\ln(2)} \, m_{\cV_{\chi}}(\cF_{c_\ell}) \leq \lim_{T \rightarrow +\infty} \frac{\# \left (k \cL_{\chi}  \cap \cE_{T} \right )}{\ln(T)}\leq \frac{\varkappa_2}{\ln(2)} \, m_{\cV_{\chi}}(\cF_{c_\ell^{-1}}).
\end{equation}
Since the function $c \mapsto m_{\cV_{\chi}}(\cF_c)$ is continuous, we have $m_{\cV_{\chi}}(\cF_{c_{\ell}^{\pm 1}}) \rightarrow m_{\cV_{\chi}}(\cF_1)$ as $\ell \rightarrow + \infty$. Let $\varkappa > 0$ be as in \eqref{eq:Varkappa-Constant}. By Lemma \ref{lem:Reduction} and letting $\ell$ go to infinity in \eqref{eq:Sandwich_lem:Reduction2}, for every $k \in K_\infty^{-1}$, we have, as $T \rightarrow + \infty$,
\begin{equation*} 
\cN_{\beta_\chi} (kx_0, T) = [K \cap L : K \cap P]^{-1} \, \# \left (k^{-1} \cL_{\chi}  \cap \cE_{T} \right ) \, \sim \, \varkappa \, \ln(T),
\end{equation*}
Noting that $X = K/(K \cap P)$ and that, for every $k \in K$ and $p \in K \cap P$, we have $\cN_{\beta_\chi} (k p x_0, T) = \cN_{\beta_\chi} (k x_0, T)$, implies that, for $\sigma_X$-almost every $x \in X$, as $T \rightarrow + \infty$,
\begin{equation*} 
\cN_{\beta_\chi} (x, T) \, \sim \, \varkappa \, \ln(T).
\end{equation*}
The proof of Theorem \ref{thm:ThmMain} is complete.
\end{proof}

\section{Proof of Theorem \ref{thm:ThmMainEquid}} \label{sec:Equidistribution-Rational}
In this section, we briefly explain how Theorem \ref{thm:ThmMainEquid} follows from the proof of Theorem \ref{thm:ThmMain2}.

\medskip
For the reader’s convenience, we recall that, for every $x \in X$ and $r > 0$, we denote by $B_X(x,r)$ the open ball in $X$ with center $x$ and radius $r$, and define the counting function 
\[
\cN_{\chi}(x,r,T) =\#\left\{v \in \bX(\Q) : v \in B_X(x,r), \, 1 \leq H_\chi(v) < T\right\}.
\]
Let $\psi_r : \R_+ \rightarrow \R_+$ be the constant function defined by $\psi_r(t) = r$. Then, in the notation of Theorems \ref{thm:ThmMain2} and \ref{thm:ThmMain}, we remark that, for every $x \in T$ and $T \geq 1$, we have $\cN_{\chi}(x,r,T) = \cN_{\psi_r}(x,T)$.
\begin{proof} [Proof of Theorem \ref{thm:ThmMainEquid}] \label{proof:Theorem-Equidistribution} 
Let the constant $\varkappa > 0$ be given as in Equation \eqref{eq:Varkappa-Constant}. Fix $x \in X$ and $r > 0$. For every $T \geq 1$, we define the set
\begin{equation} \label{eq:Definition_E_r}
\cE(r,T) = \left \{ \bm{v} \in \cV_{\chi} : d(x_0, [\bm{v}]) < r, \, 1 \leq \|\bm{v}\| < T \right \}.
\end{equation}
By Lemma \ref{lem:Reduction}, for every $x \in X$ and $T \geq 1$, we have
\[
\cN_{\chi}(x, r, T) = [K \cap P : K \cap L]^{-1} \, \# \left (k_x^{-1} \cL_{\chi}  \cap \cE(r,T) \right ).
\]
For $T \geq 1$, we define 
\[
\cF(r,T) = \left \{ \bm{v} \in \cV_{\chi}: d(x_0, [\bm{v}]) < r, \, \max \{1 ,{T}/2\} \leq \|\bm{v}\| < {T} \right \}.
\]
We decompose $\cE(r,T)$ dyadically:
\[
\cE(r,T) = \bigsqcup_{j \geq 0} \, \cF(r,T_j) \quad \text{ where } \,  T_j = T /2^{j}, \, j \geq 0.
\]
\medskip
For any $\delta \in (0,1)$, we define the sets
\[
\underline{\cF}(r,T,\delta) = \left\{ \bm{v} \in \cV_{\chi}: d(x_0, [\bm{v}]) < (1 + \delta)^{- 1} \, r, \, (1 + \delta)\tfrac{T}{2} < \|\bm{v} \| < (1 + \delta)^{-1} T \right\}
\]
and
\[
\overline{\cF}(r, T,\delta) = \left\{ \bm{v} \in \cV_{\chi}: d(x_0, [\bm{v}]) < (1 + \delta) \, r, \, (1 + \delta)^{-1}\tfrac{T}{2} < \|\bm{v} \| < (1 + \delta) T \right\}.
\]
Then, following the notation of Lemma \ref{lem:WellRounded}, we simply put (without applying an element of the subgroup $A$)
\[
\cB(r, T) = \cF(r, T), \quad \underline{\cB}(r, T,\delta) = \underline{\cF}(r, T,\delta) \quad \text{and} \quad \overline{\cB}(r,T,\delta) = \overline{\cF}(r,T,\delta).
\]
Then Lemma \ref{lem:WellRounded} and Proposition \ref{prop:WellRounded} still apply, and hence the family $(\cB(r,T))_{T \geq 1}$ is well-rounded in the sense of Definition \ref{def:Well-Rounded}. In the proof of Proposition \ref{prop:Counting}, we used the ergodicity of the action of $A$ on $G / \G$ to show that for almost every $x \in X$, we have
\[
\lambda_1(  \Ad(a_{y_T} k_x^{-1}) \kg_\G ) = T^{o(1)} \quad \text{ as }T \to + \infty.
\]
Since here $y_T = 1$, and hence $a_{y_T} = 1$, we have, for every $x \in X$, that 
\[
\lambda_1(  \Ad(a_{y_T} k_x^{-1}) \kg_\G ) \asymp 1,
\]
and hence also $\mathrm{ht}(a_{y_T} k_x^{-1} \G) \asymp 1$. Hence, by Lemma \ref{lem:Asymptotic}, there exist constants $\varkappa_2 > 0$, $\varepsilon_2 > 0$, and $T_0 > 1$ such that for all $T > T_0$,
\[
\# \left (k_x^{-1}\cL_{\chi} \cap \cB(r,T) \right ) = \varkappa_2 \, m_{\cV_{\chi}}(\cB(r,T)) \left (1 + O( T^{-\varepsilon_2}) \right ).
\]
The proof now proceeds similarly (with $\tau = 0$) as that of Proposition \ref{prop:Counting}. Thus, there exist constants $\varkappa_1 > 0$ and $\varepsilon_1 > 0$, such that
\[
\# \left (k_x^{-1}\cL_{\chi} \cap \cE(r,T) \right ) = \varkappa_1 \, m_{\cV_{\chi}}(\cE(r,T)) \left (1 + O( T^{-\varepsilon_1}) \right ).
\]
By a similar volume computation as that in the proof of Theorem \ref{thm:ThmMain2}, we have 
\begin{align*}
m_{\cV_{\chi}}(\cE(r,T)) &= \beta_\chi \,  \int_{1}^{T} \sigma_X(B_X(r)) y^{\beta_\chi d} \, \omega_0 \, \frac{d y}{y} \\ \notag
&= \frac{1}{d} \, \omega_0 \, \sigma_X(B_X(r)) \ T^{\beta_{\chi} d} + O(1),
\end{align*} 
as required. This establishes the proof of Theorem \ref{thm:ThmMainEquid} with $\varepsilon' = \varepsilon_1$ and
\begin{equation} \label{eq:Varkappa-Constant-Equid}
\varkappa' = [K \cap P:K \cap L]^{-1} \, \varkappa_1 \, \frac{1}{d} \, \omega_0.
\end{equation}
\end{proof}

\section{Applications}
Let us briefly outline how our applications follow from Theorems \ref{thm:ThmMain2} and \ref{thm:ThmMain}.

\subsection{Projective quadrics} \label{sec:Quadric}
Let $n \geq 1$ be a positive integer, and let $Q : \R^{n+2} \rightarrow \R$ be a non-degenerate rational quadratic form in ${n+2}$ variables. 
We write $\bX_{Q} = \bG / \bP$, where $\bG = \SO_Q$ is the special orthogonal group associated to the quadratic form $Q$, and $\bP \leq \bG$ is the parabolic $\Q$-subgroup stabilizing a rational isotropic line in the standard representation. If the quadratic form $Q$ is not conjugate over $\Q$ to the exceptional quadratic form $Q_0(\bm{x}) = x_1 x_4 - x_2 x_3$, then the group $\bG$ is $\Q$-simple and the parabolic $\Q$-subgroup $\bP$ is maximal, that is, $\mathrm{rank}_{\Q} \, \bP = \mathrm{rank}_{\Q} \, \bG - 1$. On the other hand, if $Q$ is conjugate over $\Q$ to $Q_0$, we have an isomorphism
\[
\bG \simeq \SO(2,2) \simeq \SO(2,1) \times \mathrm{SO}(2,1),
\]
and $\mathrm{rank}_{\Q} \, \bP = \mathrm{rank}_{\Q} \, \bG - 2$. The space of real points of $\bX_{Q}$ are identified with the $n$-dimensional projective rational quadric hypersurface given as the set of zeros in $\mathbb{P}(\R^{n+2})$ of the quadratic form $Q$:
\begin{equation}
X_Q = [Q^{-1}(0)] = \big\{ x \in \mathbb{P}(\R^{n+2}) : x = [\bm{x}]\text{ with }Q(\bm{x}) = 0 \big\}.
\end{equation}
The distance $d(\cdot,\cdot)$ and the height function $M$ are obtained by restriction of the usual distance and height function on $\mathbb{P}(\R^{n+2})$, respectively. Let $K$ be a maximal compact subgroup of the special orthogonal group $\SO_Q(\R)$ associated to $Q$ and let $\sigma_{Q}$ be the $K$-invariant probability measure on $X_Q$. Furthermore, one assumes that $X_Q$ contains a rational point; by stereographic projection, this implies in fact that $\bX_Q(\Q)$ is dense in $X_Q$.

\medskip
As mentioned in Remark~\ref{remark:Exceptional}, let us briefly explain why our methods, both for counting below and at the Diophantine exponent, do not apply to the exceptional quadric hypersurface $X_0$. Indeed, Lemmas~\ref{lem:Asymptotic} and~\ref{lem:Equi2} rely on the assumption that the parabolic subgroup defining $X$ is maximal, whereas the parabolic subgroup defining $X_0$ is \emph{not} maximal. This maximality is essential for the existence of a $G$-invariant measure on the cone $\cV_{\chi}$, for guaranteeing that $L \cap \Gamma$ is a lattice in $L$, and for the validity of the mean value formula~\eqref{eq:Mean-Value-Formula}, which was also used in \eqref{eq:First_Moment}.

\subsection{Grassmann varieties}
Let $1 \leq \ell < n$ be positive integers, and let us write $\bX_{\ell} = \bG / \bP$, where $\bG = \SL_n$ and $\bP \leq \bG$ is the parabolic $\Q$-subgroup stabilizing a rational line spanned by a pure tensor in the $\ell$-th exterior power of the standard representation of $\bG$. Then the parabolic $\bP$ is maximal and $X_{\ell}$, viewed as a subvariety of $\mathbb{P}(\bigwedge^\ell \R^n)$, is the Grassmann variety $\mathrm{Gr}_{\ell,n}(\R)$ of $\ell$-dimensional subspaces of $\R^n$. This is in accordance with Schmidt's paper \cite{Schmidt67}, where he used the Pl\"ucker embedding to define the height $H(v)$ of a rational subspace $v$ of $\R^n$. The distance used on $X_{\ell}$ is the usual Riemannian distance and we equip $X_{\ell}$ with the unique probability measure $\sigma_{\ell}$ invariant under $K=\SO_n(\R)$. Let $\bT$ be the subgroup of $\bG$ consisting of all diagonal matrices. Then $\bT$ is a maximal $\Q$-split $\Q$-torus. Let $\bP_0$ be the Borel subgroup of $\bG$ consisting of all upper-triangular matrices. Let $\Phi(\bG, \bT)$ be the associated root system with ordering induced by $\bP_0$, $\Delta$ the set of simple roots, and $\{\lambda_{\alpha}\}_{\alpha \in \Delta}$ the set of fundamental $\Q$-weights. Let $\alpha \in \Delta$ be the simple root such that $\bP = \bP_{\Delta \smallsetminus \{\alpha\}}$ is the standard parabolic $\Q$-subgroup of $\bG$ corresponding to the subset $\Delta \smallsetminus \{\alpha\}$ of simple roots. The representation of $\bG$ given by the $\ell$-th exterior power is the unique strongly rational representation of $\bG$ associated to the choice of dominant $\Q$-weight $\chi$ given by the fundamental $\Q$-weight $\lambda_{\alpha}$.

\bibliographystyle{plain}

\end{document}